\RequirePackage{fix-cm}
\documentclass[smallcondensed]{svjour3}     
\smartqed  
\usepackage{graphicx}
\usepackage{amsmath,amssymb, amsfonts,bm,tikz,url}
\usepackage{xcolor}
\usetikzlibrary{hobby}

\usepackage{anyfontsize}
\usepackage{enumerate}
\usepackage{esint}
\usepackage{etoolbox}
\usepackage{isomath}
\usepackage{mathrsfs}
\usepackage{mathtools}
\usepackage{multicol}
\usepackage{pgfplots}
\usepackage{scalerel}
\usepackage{stackengine}
\usepackage{tabulary}
\usepackage{tikz}

\makeatletter
\newcommand*\bdot{\mathpalette\bdot@{.65}}
\newcommand*\bdot@[2]{\mathbin{\vcenter{\hbox{\scalebox{#2}{$\m@th#1\bullet$}}}}}
\makeatother

\makeatletter
\newcommand*\bddot{\mathpalette\bddot@{.65}}
\newcommand*\bddot@[2]{\mathbin{\vcenter{\hbox{\scalebox{#2}
    {$\m@th#1\smash{{}_{\bullet}^{\bullet}}$}}}}}
\makeatother

\newcommand{\circled}[2][]{%
  \tikz[baseline=(char.base)]{%
    \node[shape = circle, draw, inner sep = .5pt]
    (char) {\phantom{\ifblank{#1}{#2}{#1}}};%
    \node at (char.center) {\makebox[0pt][c]{#2}};}}
\robustify{\circled}

\newcommand{\ve}[1]{\vectorsym{#1}}

\newcommand{\vf}{\ve{f}}

\newcommand{\vu}{\ve{u}}

\newcommand{\vx}{\ve{x}}
\newcommand{\vy}{\ve{y}}

\stackMath
\newcommand\reallywidecheck[1]{%
\savestack{\tmpbox}{\stretchto{%
  \scaleto{%
    \scalerel*[\widthof{\ensuremath{#1}}]{\kern-.6pt\bigwedge\kern-.6pt}%
    {\rule[-\textheight/2]{1ex}{\textheight}}
  }{\textheight}%
}{0.5ex}}%
\stackon[1pt]{#1}{\scalebox{-1}{\tmpbox}}%
}
\usepackage[utf8]{inputenc}
\usetikzlibrary{decorations.pathreplacing,calligraphy}
\usepackage{bm}
\usepackage{mathrsfs}
\usepackage{amsmath}
\usepackage{amsfonts}
\usepackage{color}
\usepackage{setspace}
\usepackage{float}
\usepackage{url}
\usepackage{multicol}
\usepackage{subfigure}
\usepackage{bm}
\usepackage{gensymb}

\usepackage{enumitem}
\usepackage{pifont}
\usepackage[top=1in,bottom=1in,left=1in,right=1in]{geometry}
\usepackage{tcolorbox}
\tcbuselibrary{breakable}
\newcommand{\real}{\mathbb{R}}

\newcommand{\mcL}{\mathcal{L}}

\usepackage{xcolor}

\def\omg{{\Omega}}

\def\omgb{\mathcal{B}\Omega}
\def\omgbb{\mathcal{B}\mathcal{B}\Omega}

\def \phib{{\boldsymbol \phi}}
\def \psib{{\boldsymbol \psi}}

\def \fb{\vf}

\def \ub{\vu}

\def \xb{\vx}

\def \zb{\mathbf{z}}

\def \yb{\vy}

\newcommand{\verti}[1]{{\left\vert #1
    \right\vert}}
\newcommand{\vertii}[1]{{\left\vert\left\vert #1
    \right\vert\right\vert}}

\begin{document}

\title{OBMeshfree: An optimization-based meshfree solver for nonlocal diffusion and peridynamics models
\thanks{Authors were supported by by the National Science Foundation under award DMS-1753031 and the AFOSR grant FA9550-22-1-0197. Portions of this research were conducted on Lehigh University's Research Computing infrastructure partially supported by NSF Award 2019035.}}

\author{Yiming Fan \and Huaiqian You \and Yue Yu}

\institute{Y. Fan \at
Department of Mathematics, Lehigh University, Bethlehem, PA\\
\email{yif319@lehigh.edu}
\and
H. You \at
Department of Mathematics, Lehigh University, Bethlehem, PA\\
\email{huy316@lehigh.edu}
\and
Y. Yu \at
Department of Mathematics, Lehigh University, Bethlehem, PA\\
\email{yuy214@lehigh.edu}
}

\date{Received: date / Accepted: date - Draft, ***, 2022}

\maketitle

\begin{abstract}
We present OBMeshfree, an Optimization-Based Meshfree solver for compactly supported nonlocal integro-differential equations (IDEs) that can describe material heterogeneity and brittle fractures. OBMeshfree is developed based on a quadrature rule calculated via an equality constrained least square problem to reproduce exact integrals for polynomials. As such, a meshfree discretization method is obtained, whose solution possesses the asymptotically compatible convergence to the corresponding local solution. Moreover, when fracture occurs, this meshfree formulation automatically provides a sharp representation of the fracture surface by breaking bonds, avoiding the loss of mass. As numerical examples, we consider the problem of modeling both homogeneous and heterogeneous materials with nonlocal diffusion and peridynamics models. Convergences to the analytical nonlocal solution and to the local theory are demonstrated. Finally, we verify the applicability of the approach to realistic problems by reproducing high-velocity impact results from the Kalthoff–Winkler experiments. Discussions on possible immediate extensions of the code to other nonlocal diffusion and peridynamics problems are provided. OBMeshfree is freely available on GitHub \cite{OBMeshfree}.

\end{abstract}

\tableofcontents

\section{Introduction}

Nonlocal models such as nonlocal diffusion and peridynamics provide a description of governing laws in terms of integral operators rather than classical differential operators \cite{bobaru2016handbook,du2011mathematical,emmrich2007analysis,parks2008implementing,seleson2009peridynamics,silling_2000,zimmermann2005continuum}. They can describe phenomena not well represented by classical partial differential equations (PDEs), especially on problems characterized by long-range interactions and discontinuities \cite{bazant2002nonlocal,bobaru2016handbook,du2013nonlocal}. As a result, applications of nonlocal models span a large spectrum of scientific and engineering problems, including subsurface transport \cite{benson2000application,katiyar2020general,katiyar2014peridynamic,schumer2003multiscaling,schumer2001eulerian}, phase transitions \cite{bates1999integrodifferential,chen2000nonlocal,dayal2006kinetics}, image processing \cite{buades2010image,d2021bilevel,gilboa2007nonlocal,lou2010image}, multiscale and homogenized systems \cite{alali2012multiscale,askari2008peridynamics,du2020multiscale,seleson2010peridynamic,silling2021propagation,silling2022peridynamic,you2022data}, turbulence \cite{bakunin2008turbulence,schekochihin2008mhd}, stochastic processes \cite{d2017nonlocal,meerschaert2019stochastic,metzler2000random}, and fracture mechanics \cite{ha2011characteristics,silling2001kalthoff,silling_2000}.

In this work, we consider nonlocal models that are characterized by a general heterogeneous nonlocal operators of the form:
\begin{equation*}
\mcL_{\delta}[u](\xb,t):=\int_{B_\delta(\xb)}K(\xb,\yb)(u(\yb,t)-u(\xb,t))d\yb,  \quad \text{with }K(\xb,\yb):=\varpi(\xb,\yb)C(\xb,\yb)\gamma(\xb,\yb).    
\end{equation*}
where $u$ is the solution we seek, $\gamma(\xb,\yb):=\gamma_\delta(\verti{\xb-\yb})$ is a (possibly singular) radial kernel which for fixed $\xb$ is supported on the ball of radius $\delta$, $B_\delta(\xb)$, and the two point functions $\varpi(\xb,\yb)$ and $C(\xb,\yb)$ allow for the description of material heterogeneity while preserving the physical consistency. $\delta$ defines the extant of nonlocal interactions, which is also referred to as a horizon. This integral form allows for the description of long-range interactions and reduces the regularity requirements on problem solutions, and hence enhances the accuracy of their modeling representations by generalizing the space of admissible solutions, which can feature discontinuities. Another important feature of such models is that when the classical continuum models still apply and with proper definitions of $\varpi(\xb,\yb)$ and $C(\xb,\yb)$, these nonlocal models can revert back to classical continuum models with heterogeneous material properties, as $\delta \rightarrow 0$. 

When discretizing the nonlocal models, it is desired to  preserve the corresponding local limit under refinement $h \rightarrow 0$, since analyzing consistency in the limit to the local solution provides a mathematically unambiguous means to understand accuracy and physical compatibility. Such a property is termed asymptotically compatible (AC) \cite{tian2014asymptotically}\footnote{For nonlocal models one often refines both $\delta$ and $h$ at the same rate under so-called $\delta$-convergence \cite{bobaru2009convergence}, to allow scalable implementations. Although in the literature a scheme is termed AC if it recovers the solution whenever $\delta,h\rightarrow 0$, in this work we adopt a practical setting and only require the $\delta$-convergence for AC.}. In recent years, there has been significant work in recent years toward establishing such discretizations \cite{d2020numerical,du2016local,hillman2020generalized,leng2019asymptotically,pasetto2018reproducing,seleson2016convergence,tao2017nonlocal,tian2014asymptotically,trask2019asymptotically,You_2019,you2020asymptotically}. Broadly, strategies either involve adopting traditional finite element shape functions and carefully performing geometric calculations to integrate over relevant horizon/element subdomains, or adopt a strong-form meshfree discretization where particles are associated with abstract measure. The former is more amenable to mathematical analysis due to a better variational setting, while the latter is simple to implement and generally faster \cite{bessa2014meshfree,silling2005meshfree}.



In this paper we focus on the meshfree approach and approximate the heterogeneous nonlocal operator as:
$$\mcL_{\delta}[u](\xb_i,t)\approx \mcL_{\delta,h}[u](\xb_i,t):=\sum_j K(\xb_i,\xb_j)(u(\xb_j,t)-u(\xb_i,t))\omega_{ij},$$
where the quadrature weights $\omega_{ij}$ are associated with a local neighborhood of particles for each discretization point $\xb_i$, generated by local optimizations to make the approximation rule exact for certain classes of functions. By defining the averaged material property field $K(\xb,\yb)$ as an analog of a series of two springs connecting the two points, our recent work  \cite{d2022prescription,fan2021asymptotically,fan2022meshfree,foss2022convergence,trask2019asymptotically,yu2021asymptotically} has provided theoretical analysis and numerical verifications on the AC property on this optimization-based quadrature rule for  heterogeneous  materials \cite{fan2021asymptotically,fan2022meshfree}. In this paper, we will provide an open-source meshfree solver and demonstrations of its convergence properties on various examples. To achieve a convergent simulation, the AC property to the local limit is only one important ingredient. Besides the consistency to the local limit, two additional features are desired in our nonlocal problem solver. Firstly, in peridynamics problems, one of its main appeals is to handle fracture problems, where free surfaces are associated with the evolution of a fracture surface. To achieve numerical consistency for problems involving fracture, one must also consider the interplay between consistency of quadrature for discrete operators and the imposition of traction loads as fracture surfaces open up and evolve \cite{lipton2014dynamic}. Second, in applications such as the particle systems with long-range interactions, the horizon size $\delta$ should be seem as a physical value and there is possibly no corresponding local limit. To preserve the physical consistency in such a scenario, the numerical convergence to the correct nonlocal limit when $h\rightarrow 0$ would be desired in the nonlocal problem solver.

Our goal is to demonstrate a comprehensive treatment of nonlocal quadrature rule, material heterogeneity, and evolving free surfaces, which is able to achieve numerical consistency to both local and nonlocal limits and capture material fracture. In particular, when no fracture occurs and the analytical solution is sufficiently smooth, the formulation should preserve the AC limit under $\delta$-convergence and the consistency to the nonlocal limit as $h\rightarrow 0$. Moreover, when fracture occurs, the formulation should be able to capture the material damage and the evolving fracture surfaces via bond breaking. This practically means that one is able to incorporate all of the necessary ingredients to perform non-trivial simulations of fracture mechanics while maintaining a scalable implementation and guaranteeing convergence. To achieve these properties, our development has two steps. First, to handle material heterogeneity and free surfaces in such a way that one preserves a limit to the relevant local problem, a unified mathematical formulation is introduced. Then, to preserving the AC limit under $\delta$-convergence and the consistency to the nonlocal limit as $h\rightarrow 0$, an optimization-based quadrature rule is employed. As a result, our method provides an efficient discretization with rigorous underpinnings for a class of nonlocal models featuring material heterogeneity and evolving fractures.


We remark that the paper is organized to first establish the mathematical formulations and provide a brief summary of theoretical underpinnings of the approach, while the second half focuses on a user manual for the code \cite{OBMeshfree} together with demonstrations on several exemplar applications. Practitioner with more applied interests may skip the first part without issue. The paper is organized as follows. 
We firstly recall the heterogeneous nonlocal diffusion and bond-based peridynamics problems in Section \ref{sec:theory}, and provide a unified mathematical formulation for handling material heterogeneity and fracture. In Section \ref{sec:disc}, the optimization-based quadrature rule is elaborated as a unified numerical approach for heterogeneous nonlocal diffusion and bond-based peridynamics, together with the treatment of material fracture. We also provide a summary of the rigorous convergence analysis for the optimization-based quadrature rule, verifying its consistency to the local limit when $h,\delta\rightarrow 0$, and to the nonlocal limit when $h\rightarrow 0$. Then, in Section \ref{sec:manual} the main structure of OBMeshfree code is shown and each code component is discussed in detail. In Section \ref{sec:test}, we demonstrate four examples as verifications and validations of the code, including three problems with manufactured solutions and one example on reproducing high-velocity impact results from the Kalthoff–Winkler experiment as an engineering-oriented problem. 
Section \ref{sec:conclusion} summarizes our results and discusses future research.

\section{Nonlocal Theory}\label{sec:theory}

In this section, we introduce the notation and describe the nonlocal models that will be useful throughout the following sections.

Let $\omg\in\real^d$, $d=1,2,3$, be a bounded open domain. We are interested in solving for functions $u:\omg\rightarrow \real$ and $\ub:\omg\rightarrow \real^d$, solutions of nonlocal diffusion and nonlocal mechanics problems, respectively. Herein, $u(\xb)$ represents the concentration of a diffusive quantity in the nonlocal diffusion problem, and $\ub(\xb)$ represents the displacement field of an object in mechanics. In nonlocal settings, every point in a domain interacts with a neighborhood of points. In this work, we further assume that such neighborhood is an Euclidean ball surrounding points in the domain, i.e., $B_\delta(\xb):=\{\yb\in\real^d:|\yb-\xb|\leq\delta\}$, where $\delta $ is the horizon. This assumption has implications on the concept of boundary conditions. In particular, unless otherwise stated, the boundary conditions should no longer be prescribed on the sharp interface, $\partial\omg$, but on a collar of thickness of at least $\delta$ surrounding the domain $\omg$, that we refer to as the nonlocal volumetric boundary domain (or simply nonlocal boundary),
$$\omgb:=\left\{ \bm x\notin\omg|\text{dist}(\bm x, \partial\omg)  <\delta\right\}.$$
This set consists of all points outside the domain that interact with points inside the domain. To define the nonlocal problems with general mixed boundary conditions, we further decompose the sharp interface $\partial\Omega$ into two parts: $\partial\Omega=\partial\Omega_D\bigcup \partial\Omega_N$, where $(\partial\Omega_D)^o\bigcap (\partial\Omega_N)^o=\emptyset$. To apply the nonlocal Dirichlet-type boundary condition, we assume that $\ub(\xb)=\ub_D(\xb)$ are provided in a layer with non-zero volume outside $\Omega$, while the free surface boundary condition is applied on the sharp interface $\partial\Omega_N$. To define a Dirichlet-type constraint, we denote
\begin{align*}
\omgb_D:=\{\xb\notin\Omega|\text{dist}(\xb,\partial\Omega_D)<\delta\},\,
\end{align*}
and assume that the value of $\ub$ is given on $\omgb_D$. For notation simplicity, we denote $\omg_D:=\omg\bigcup\omgb_D$.

In this paper and the code \cite{OBMeshfree}, we focus on 2D problems ($d = 2$) and provide demonstrations with both static and dynamic examples, although the method is also applicable to 3D problems.

\subsection{Nonlocal Diffusion Models}\label{sec:theory_diff}

Nonlocal diffusion models have been employed in many applications \cite{bakunin2008turbulence,bobaru2010peridynamic}, and they are capable to describe the underlying phenomena when the classical Fick's first law or standard Brownian motion fails \cite{bucur2016nonlocal,du2012analysis,metzler2004restaurant,neuman2009perspective}. Specifically, given a loading field $f$, the time-dependent nonlocal diffusion equation can be given as:
\begin{equation}\label{eqn:diffusion}
    \dfrac{\partial u}{\partial t}(\xb,t)-\mcL_{D\delta}[u](\xb,t) =f(\xb,t),
\end{equation}
where the diffusion operator on scalar function $u:\real^d\rightarrow \real$ is defined as 
\begin{equation}\label{eq:diffusion_comp}
\mcL_{D\delta}[u](\xb,t):=2\int_{B_\delta(\xb)}A(\xb,\yb)\gamma(\xb,\yb)(u(\yb,t)-u(\xb,t))d\yb=f(\xb,t), \quad \xb\in\omg.
\end{equation}
Here, $A(\xb,\yb)\in[r,R]$, $r,R>0$, is a uniformly bounded and continuous two-point function, describing the heterogeneous diffusion property. $\gamma$ denotes a properly scaled kernel function which is assumed to satisfy:
\begin{equation}\label{eqn:require_ga}
\gamma(\xb,\yb)=\gamma_\delta(|\xb-\yb|)=\frac{1}{\delta^{d+2}}\gamma_1\left(\frac{|\xb-\yb|}{\delta}\right)=\frac{D_0}{\delta^{d+2-s}|\xb-\yb|^s},
\end{equation}
where $\gamma_1$ is a nonnegative and nonincreasing function, and there exists a positive constant $\zeta<1$ satisfying $B_\zeta(\bm{0})\subset\text{supp}(\gamma_1)\subset B_1(\bm{0})$ and $\int_{B_1(\bm{0})}\gamma_1(|\zb|)|\zb|^2d\zb=d$. Then, \eqref{eqn:diffusion} can be seen as a nonlocal analogue to the local diffusion equation\footnote{Herein, we apply the assumptions in \eqref{eqn:require_ga} so as to guarantee the compatibility of the nonlocal model and its local limit, in the limit of vanishing nonlocality, i.e., $\delta\rightarrow 0$. However, we point out that the optimization-based quadrature rule as well as the OBMeshfree package can be readily applied to more general kernels, such as the data-driven kernel developed in \cite{you2022data}.}. In particular, when taking the local diffusion parameter field $a(\xb):=A(\xb,\xb)$ and consistent Dirichlet-type boundary conditions, it was shown in \cite{fan2021asymptotically} that \eqref{eqn:diffusion} is well-posed and converges to the classical diffusion equation
\begin{equation}\label{eqn:diffusion_local}
    \dfrac{\partial u}{\partial t}(\xb,t)-\mcL_{D}[u](\xb,t) =f(\xb,t),\quad \mcL_{D}[u](\xb,t):=\nabla\cdot(a(\xb)\nabla u(\xb,t)),
\end{equation}
as $\delta\rightarrow 0$. Therefore, in examples where only the local diffusion coefficient field $a(\xb)$ is provided, one can take the nonlocal diffusion coefficient as the harmonic mean of the local diffusion coefficient:
\begin{equation}
A(\xb,\yb)=2\left(\dfrac{1}{a(\xb)}+\dfrac{1}{a(\yb)}\right)^{-1}    
\end{equation}
For further details, we refer interested readers to \cite{fan2021asymptotically}.

Here, we consider nonlocal diffusion problems with Dirichlet-type boundary conditions without loss of generality\footnote{Neumann and Robin-type boundary conditions can be implemented with the quadrature rule provided by OBMeshfree, following \cite{d2020physically,You_2019,you2020asymptotically}}. That means, $\omgb=\omgb_D$ and a volume constraint $u_D:\omgb\times[0,T]\rightarrow \real$ is provided. Then, the nonlocal counterpart of a Dirichlet boundary condition for PDEs is applied as a volume constraint: $u(\xb,t)=u_D(\xb,t)$ for $(\xb,t)\in \omgb\times[0,T]$. Although in Section \ref{sec:test} we only provide numerical verification for the convergence of numerical solutions in static cases, sample codes for both static and time-dependent nonlocal diffusion problems are provided in \cite{OBMeshfree}.

\subsection{Peridynamics Models}\label{sec:theory_peri}

The peridynamic theory provides a nonlocal mechanics model, which has been applied for material failure and damage simulation \cite{bobaru2016handbook,diehl2022comparative,diehl2019review,silling_2010,silling_2000,silling2003dynamic} and provided robust modeling capabilities for analysis of complex crack propagation phenomena, such as crack branching \cite{bobaru2015cracks,ha2010studies,ha2011characteristics,yu2021asymptotically}, bridging, deflection and trapping \cite{fan2022meshfree,prakash2022investigation}.

Consider a body occupying the domain $\omg$, the general peridynamic equation of motion for a point $\xb\in\omg$ and time $t\in[0,T]$ is
$$\rho\dfrac{\partial^2\ub}{\partial t^2}(\xb,t)-\mcL_{P\delta}[\ub](\xb,t)=\fb(\xb,t),$$
where $\mcL_{P\delta}$ is a nonlocal operator representing the peridynamic internal force density, $\rho$ is the mass density, and $\fb$ is a prescribed body force density. As for the nonlocal diffusion problems, the nonlocal interactions in $\mcL_{P\delta}$ are also restricted to the nonlocal neighborhood, $B_\delta(\xb)$, characterized by the horizon size $\delta$. In this work, we focus on the bond-based peridynamic solid model \cite{bobaru2016handbook,trask2019asymptotically,Yu2018paper}, and take the peridynamics operator as:
\begin{equation}\label{eq:nonlocElasticity_comp}
\begin{aligned}
    \mcL_{P\delta}[\ub](\xb,t):=&
  c\int_{B_\delta (\xb)} \kappa(\xb,\yb) \gamma(\left|\yb-\xb\right|)\frac{\left(\yb-\xb\right)\otimes\left(\yb-\xb\right)}{\left|\yb-\xb\right|^2}  \left(\ub(\yb,t) - \ub(\xb,t) \right) d\yb.
  \end{aligned}
\end{equation}
Here, $\gamma$ is the kernel function as defined in \eqref{eqn:require_ga}, and the two-point functions $\kappa(\xb,\yb)$ denote the (averaged) bulk modulus property\footnote{In 2D problems, the bond-based peridynamics model is restricted with Poisson ratio $\nu=0.25$, and in 3D problems one has a fixed Poisson ratio $\nu=1/3$. Hence, the bulk modulus, $\kappa$, and the Young's modulus, $E$, can be converted from each other following  $\kappa=\frac{E}{3(1-2\nu)}$.}. To recover parameters for linear elasticity when the nonlocal effects vanish, one should take $c=24/5$ for $d=2$, and $c=6$ for $d=3$ (see, e.g., \cite{emmrich2007analysis,trask2019asymptotically} for further details). Similar as in the nonlocal diffusion problems, in examples where only the local bulk modulus coefficient field $\kappa(\xb)$ is provided, we take the nonlocal bulk modulus coefficient as the harmonic mean of the local coefficient \cite{mehrmashhadi2018effect,nguyen2021depth,prakash2022investigation}: 
\begin{equation}\label{eqn:hamo_k}
\kappa(\xb,\yb)=2\left(\dfrac{1}{\kappa(\xb)}+\dfrac{1}{\kappa(\yb)}\right)^{-1}.    
\end{equation}




One of the main features of peridynamics is to handle fracture problems, where damage is incorporated into the peridynamic constitutive model by allowing the bonds of solid interactions to break irreversibly. Here we employ the critical stretch criterion where breakage occurs when a bond is extended beyond some predetermined critical bond deformed length \cite{yu2021asymptotically,zhang2018state}. Then, this criterion is implemented by multiplying the pairwise force function with a history-dependent scalar boolean function. In particular, to model brittle fracture, we break the bond between two material points, $\xb$ and $\yb$, when the associated strain exceeds a critical strain criteria. 
A two-point boolean state function $\theta(\xb,\yb,t)$ is defined and updated to describe the bond breakage through the crack growing 
\begin{align}\label{eqn:theta}
    \theta(\xb,\yb,t) &= \begin{cases}
    1, \quad \text{if }s(\xb,\yb,\tau)\leq s_0(\xb,\yb),\;\forall\tau\leq t, \text{ and }\yb\in B_\delta(\xb)\bigcap\Omega_D,\\
    0, \quad \text{otherwise}, \\
    \end{cases}
\end{align}
where the 
associated strain $s$ and the critical bond stretch $s_0$ related to material parameters are defined as \cite{zhang2018state}:
\begin{align}
 \nonumber&s(\xb,\yb,t): = \frac{||\ub(\yb,t)-\ub(\xb,t) + \yb-\xb||-||\yb - \xb||}{||\yb - \xb||},\\
&s_0(\xb,\yb) := \begin{cases}
    \sqrt{\frac{\pi G(\xb,\yb)}{3\kappa(\xb,\yb)\delta}}, \quad \text{2D},\\
    \sqrt{\frac{5G(\xb,\yb)}{9\kappa(\xb,\yb)\delta}}, \quad \text{3D}, \\
    \end{cases}\label{eq:criteria}
\end{align}
where $\kappa$ is the nonlocal bulk modulus coefficient. 
In \eqref{eq:criteria} $G(\xb,\yb)$ is the averaged fracture energy defined via the arithmetic mean:
\begin{equation}\label{def_fract_energy}
   G(\xb,\yb)= \frac{1}{2}(G(\xb)+G(\yb)).
\end{equation}
Here, the averaged material properties definition in \eqref{eqn:hamo_k} and the averaged fracture energy definition in \eqref{def_fract_energy} are inspired by seeing the interaction between $\xb$ and $\yb$ as an analog of a series of two springs connecting the two points. Following a similar argument as in \cite{fan2022meshfree}, one can show that when no fracture occurs and the local modulus field satisfies $\kappa\in C(\overline{\omg\bigcup\omgb})$, it is guaranteed that the nonlocal solution of \eqref{eq:nonlocElasticity_comp} converges to the solution of a linear elastic model as $\delta\rightarrow 0$, hence it preserves the correct local limit.

In summary, with proper initial conditions and Dirichlet-type boundary condition in $\omgb_D$, we obtain a unified mathematical formulation for dynamic bond-based peridynamics for $\xb\in\omg$:
\begin{align}
    \rho\dfrac{\partial^2\ub}{\partial t^2}(\xb,t)
-  c\int_{B_\delta (\xb)} \theta(\xb,\yb,t)\kappa(\xb,\yb)
     \gamma(\left|\yb-\xb\right|)\frac{\left(\yb-\xb\right)\otimes\left(\yb-\xb\right)}{\left|\yb-\xb\right|^2} 
      \left(\ub(\yb,t) - \ub(\xb,t) \right) d\yb= \fb(\xb,t).\label{eq:newform2}
\end{align}
This formulation naturally handles both material heterogeneity and evolving fracture as the free surface boundary conditions on $\partial\omg_N$. For further discussions on the free surface boundary conditions and more general Neumann-type boundary conditions in peridynamics, we refer interested readers to \cite{fan2022meshfree,yu2021asymptotically}.

\section{Optimization-Based Quadrature Rules}\label{sec:disc}


In this section, we elaborate the strong-form particle discretization of the nonlocal diffusion and peridynamics models introduced above. To obtain a unified formulation for both models, we rewrite \eqref{eq:diffusion_comp} and \eqref{eq:nonlocElasticity_comp} as a general heterogeneous nonlocal operator of the form:
\begin{equation}\label{eqn:general_nonlocal}
    \mcL_{\delta}[u](\xb,t):=\int_{B_\delta(\xb)}K(\xb,\yb)(u(\yb,t)-u(\xb,t))d\yb,  \quad \text{with }K(\xb,\yb):=\varpi(\xb,\yb)C(\xb,\yb)\gamma(\xb,\yb).
\end{equation}
Here, $\gamma$ is the nonlocal kernel satisfying \eqref{eqn:require_ga}, $\varpi(\xb,\yb)$ is a general two-point function corresponding to the (heterogeneous) material properties, and $C(\xb,\yb)$ corresponds to a two-point (tensor) function which is designed to guarantee the consistency to the desired local limit. In nonlocal diffusion problems, we have $\varpi(\xb,\yb):=A(\xb,\yb)$ and $C(\xb,\yb):=2$, which corresponds to the diffusion property. In bond-based peridynamics, $\varpi(\xb,\yb):=\kappa(\xb,\yb)$ and $C(\xb,\yb):=c\dfrac{(\yb-\xb)\otimes(\yb-\xb)}{\verti{\yb-\xb}^2}$, characterizing the average material properties and bond strengths between material points $\xb$ and $\yb$.


Denoting $u_{\delta}$ and $u_0$ as the nonlocal and local analytical solution respectively, and $u_{\delta,h}$ as the numerical solution, in OBMeshfree we focus on two types of convergence:
\begin{equation}\label{eqn:convh}
\underset{h\rightarrow 0}{\lim}\vertii{u_{\delta,h}-u_\delta}_{L^2(\omg)}=0,\quad \text{ and }\quad\underset{h,\delta\rightarrow 0}{\lim}\vertii{u_{\delta,h}-u_0}_{L^2(\omg)}=0.
\end{equation}
The first type of convergence indicates that the numerical discretization method is consistent with the nonlocal problem, while the second type shows that the nonlocal numerical solution preserves the correct local limit, or equivalently, the numerical scheme is asymptotically compatible. To maintain an easily scalable implementation, in asymptotic compatibility studies we assume $\delta$ to be chosen such that the ratio $\frac{\delta}{h}$ is bound by a constant as $\delta \rightarrow 0$, restricting ourselves to the ``$\delta$-convergence'' scenario \cite{bobaru2009convergence}. In the following sections, we will first introduce the spatial and temporal discretization methods for nonlocal diffusion and peridynamics problems with full Dirichlet-type boundary conditions. Then, we incorporate the bond-breaking mechanism and the free surface formulation, to provide a fully discretized framework for bond-based peridynamics including the damage criteria and the handling of free surfaces created by evolving fracture.


\subsection{Spatial Discretization}

Assume that the whole interaction domain, $\omg_D$, is discretized into a collection of points 
$$\chi_{h} = \{\xb_i\}_{i=1}^M \subset \Omega\bigcup\omgb,$$
we seek for numerical solutions such that $u_i\approx u(\xb_i)$. Recall the definitions \cite{wendland2004scattered} of fill distance $h_{\chi_h,\Omega} = \underset{\xb_i \in \Omega\bigcup\omgbb}{\sup}\, \underset{\xb_j \in \chi_h}{\min}||\xb_i - \xb_j||_2$ and separation distance
${q_{\chi_h} = \frac12 \underset{i \neq j}{\min} ||\xb_i - \xb_j||_2}$, for simplicity we drop subscripts and simply write $h$ and $q$. In this work, we assume that $\chi_h$ is \textit{quasi-uniform}, namely that there exists a constant $C_{q} > 0$, such that $q \leq h \leq C_q q$. 
Then, we seek to generate consistent meshfree quadrature rules of the form
\begin{equation}\label{eqn:quad}
\mcL_{\delta}[u](\xb_i,t)\approx \mcL_{\delta,h}[u](\xb_i,t):=\sum_j K(\xb_i,\xb_j)(u(\xb_j,t)-u(\xb_i,t))\omega_{ij}.
\end{equation}
Here, $\{\omega_{ij}\}_{\xb_j \in B_\delta (\xb_i)}$ is a collection of to-be-determined quadrature weights corresponding to a neighborhood of collocation point $\xb_i$, which will be constructed through an optimization-based approach to ensure consistency guarantees. Specifically, we seek quadrature weights for integrals supported on balls of the form
\begin{equation}
I[q] := \int_{B_\delta (\xb_i)} q(\xb_i,\yb) d\yb \approx I_h[q] := \sum_{\xb_j \in \chi_h\bigcap B_\delta(\xb_i)\backslash\{\xb_i\}} q(\xb_i,\xb_j)\omega_{ij}
\end{equation}
where the subscript $i$ in $\left\{\omega_{ij}\right\}$ denote that we seek a different family of quadrature weights for different subdomains $B_\delta(\xb_i)$. These weights are then generated from the following optimization problem
\begin{align}\label{eq:quadQP}
  \underset{\left\{\omega_{ij}\right\}}{\text{min}} \sum_{\xb_j \in \chi_h\bigcap B_\delta(\xb_i)\backslash\{\xb_i\}} \omega_{ij}^2\gamma_\delta(\xb_i,\xb_j) \quad
  \text{such that}, \quad
  I_h[q] = I[q] \quad \forall q \in \bm{V}_{h,\xb_i},
\end{align}
where  $\bm{V}_{h,\xb_i}=\left\{q(\yb-\xb_i)=p(\yb-\xb_i)\gamma_\delta(\xb_i,\yb)C(\xb_i,\yb)\Big|p\in\mathbb{P}_n(\mathbb{R}^d)\text{ such that } \int_{B_\delta(\xb_i)}q(\yb-\xb_i)d\yb<\infty\right\}$ denotes the space of functions which should be integrated exactly. $\mathbb{P}_n(\real^d)$ is the space of $n$-th order polynomials, and $C(\xb,\yb)$ is given as in \eqref{eqn:general_nonlocal}.

For each material point $\xb_i\in\omg\bigcap \chi_h$, we denote the total number of to-be-determined quadrature weights $\omega_{ij}$ as $M_i$. To solve for the optimization problem \eqref{eq:quadQP}, we formulate it as a saddle point problem 
\begin{equation}\label{eqn:QPsaddle}
\begin{bmatrix}
\bm{W} & \bm{B}^\intercal \\
\bm{B} & \bm{0}
\end{bmatrix}
\begin{bmatrix}
\bm{\omega} \\
\bm{\lambda}
\end{bmatrix}
 =
 \begin{bmatrix}
 \bm{0} \\
 \bm{g}
 \end{bmatrix},
\end{equation}
where 
$\bm{W}\in \real^{M_i\times M_i}$ is a diagonal matrix with the diagonal element determined by $\bm{W}_{j,j}=2\gamma_\delta(\xb_i,\xb_j)$, $\bm{\omega} \in \real^{M_i}$ are the vector of quadrature weights $\omega_{ij}$, and $\bm{\lambda}\in \real^{\text{dim}(V_{h,\xb_i})}$ are a set of Lagrange multipliers. $\bm{B}\in \real^{M_i\times \text{dim}(V_{h,\xb_i})}$ consists of the reproducing function evaluated at the quadrature points, satisfying $\bm{B}_{\alpha,j}=q_{\alpha}(\xb_j-\xb_i)$, for $q_{\alpha}\in \bm{V}_{h,\xb_i}$ and $\xb_j\in \chi_h\bigcap B_\delta(\xb_i)\backslash\{\xb_i\}$. $\bm{g}\in \real^{\text{dim}(V_{h,\xb_i})}$ consists of the integral of the reproducing functions over the ball, satisfying  $\bm{g}_\alpha=I[q_\alpha]$. By eliminating the constraints, the quadrature weights can be obtained by solving
\begin{equation}\label{eqn:linearomega}
\bm{\omega} = \bm{W}^{-1}\bm{B}^\intercal(\bm{B}\bm{W}^{-1}\bm{B}^\intercal)^{-1}\bm{g}.
\end{equation}
We note that the application of this quadrature does not require a background grid and is therefore truly meshfree. Moreover, the quadrature weights only need to be obtained once, using a list of neighbors lying within $B_\delta(\xb)$. In fact, in OBMeshfree weights are obtained as a preprocessing step, by solving a small local optimization problem requiring only the inversion of a small linear system in \eqref{eqn:linearomega} for each $\xb_i$.


Substituting the quadrature rule \eqref{eqn:quad}, the nonlocal diffusion operator $\mcL_{D\delta}$ and the bond-based peridynamics operator $\mcL_{P\delta}$ can be respectively discretized as:
\begin{align}
    (\mathcal{L}_{D\delta,h} u)_i&:=2\sum_{\xb_j \in \chi_h\bigcap B_\delta(\xb_i)\backslash\{\xb_i\}}A(\xb_i,\xb_j)\gamma_\delta(|\xb_i-\xb_j|)(u_i-u_j)\omega_{ij},\label{eqn:discreteNonlocDiff}\\
    (\mathcal{L}_{P\delta,h} \ub)_i&:=c\sum_{\xb_j \in \chi_h\bigcap B_\delta(\xb_i)\backslash\{\xb_i\}} \kappa(\xb_i,\xb_j) \gamma_\delta(|\xb_i-\xb_j|)\frac{\left(\xb_j-\xb_i\right)\otimes\left(\xb_j-\xb_i\right)}{\left|\xb_j-\xb_i\right|^2}  \left(\ub_i - \ub_j \right) \omega_{ij}.\label{eqn:discreteNonlocBond}
\end{align}
As verified in \cite{fan2021asymptotically,fan2022meshfree,foss2022convergence,trask2019asymptotically}, the above quadrature rule is able to obtain a compatible discretization that achieves both the consistency with the nonlocal analytical solution and the asymptotic compatibility to the local limit. In the following, we briefly list the theoretical analysis results provided in \cite{fan2021asymptotically,trask2019asymptotically}. For analysis, we denote 
\begin{equation}\label{eqn:require_k}
K(\xb,\yb)=\dfrac{N(\xb,\yb)}{\delta^{d+2-s}|\xb-\yb|^s},
\end{equation}
where the numerator satisfies $N(\xb,\yb)\leq C_N$ for all $\yb\in B_\delta(\xb)$.

For the asymptotic compatibility to the local limit, as shown in \cite{trask2019asymptotically}, the optimization-based quadrature rule has the following truncation error estimates:
\begin{theorem}[Truncation Estimates with fixed $h/\delta$]\label{thm:trunc_err}
Consider a fixed ratio $h/\delta$ and assume that both $N(\xb,\yb)$ and $u(\xb)$ are sufficiently smooth, i.e., $N(\xb,\yb)\in C^{n+2}(\overline{(\omg_D)^2})$ and $u(\xb)\in C^{n+2}(\overline{\omg_D})$. For any $\xb_i\in\chi_h\bigcap\omg$, the quadrature weights obtained from \eqref{eq:quadQP} with $n>d+s-3$ would satisfy the following pointwise error estimate, with a constant $C>0$ independent of $\delta$ and $\xb_i$:
$$\left| \int_{B_\delta(\xb_i)} K(\xb_i,\yb)u(\yb)d\yb
-\sum_{\xb_j \in \chi_h\bigcap B_\delta(\xb_i)\backslash\{\xb_i\}} K_{ij}u_j \omega_{ij}\right|<C\delta^{n-1}$$
\end{theorem}

For the convergence to the nonlocal analytical solution, following \cite{fan2021asymptotically}, we further require that $s<d$. Moreover, $\chi_h$ is assumed to be a uniform Cartesian grid:
$$\chi_h:=\{(k_{(1)}h,\cdots,k_{(d)}h)|\bm{k}=(k_{(i)},\cdots,k_{(d)})\in\mathbb{Z}^d\}\bigcap \omg_D,$$
where $h$ is the spatial grid size. Then, the optimization-based quadrature rule has the following truncation error estimates with fixed $\delta$ and vanishing $h$:
\begin{theorem}[Truncation Estimates with fixed $\delta$]\label{thm:trunc_err_nonlocal}
Consider a fixed $\delta$ and assume that both $N(\xb,\yb)$ and $u(\xb)$ are sufficiently smooth, satisfying $N(\xb,\yb)\in C^{4}(\overline{(\omg_D)^2})$ and $u(\xb)\in C^{1}(\overline{\omg_D})$. Then there exists a constant $C_{pos}<1$, such that for $h<C_{pos}\delta$, the quadrature weights obtained from \eqref{eq:quadQP} with $n=3$ would satisfy the following pointwise error estimate, with the generic constant $C>0$ independent of $h$ but may dependent on $\delta$:
$$\left| \int_{B_\delta(\xb_i)} K(\xb_i,\yb)u(\yb)d\yb
-\sum_{\xb_j \in \chi_h\bigcap B_\delta(\xb_i)\backslash\{\xb_i\}} K_{ij}u_j \omega_{ij}\right|<Ch^{\text{min}(1,d-s)}.$$
\end{theorem}

With the above truncation error estimates, the following convergence properties can be proved for static nonlocal diffusion problems, with detailed proof can be found in \cite{fan2021asymptotically}:

\begin{theorem}[Asymptotic Compatibility to the Analytical Local Diffusion Problems]\label{thm:AC}
Consider uniform Cartesian grids and $s<d$. Assume that $A(\xb,\yb)\in C^{4}(\overline{(\omg_D)^2})$, $a(\xb)\in C^{\infty}(\omg)$, $\omg_D\in C^1$, and the analytical local diffusion solution $u_0\in C^4({\overline{\omg_D}})$. When applying the boundary condition $u_D(\xb_i)=u_0(\xb_i)$ for $\xb_i\in\chi_h\bigcap\omgb_D$, there exists a $\delta_0>0$ and $C_{pos}<1$, such that for any $0<\delta<\delta_0$ and fixed ratio $h/\delta<C_{pos}$, the meshfree quadrature rule with $n=3$ is asymptotically compatible for nonlocal diffusion problems, i.e.,
\begin{equation}
\vertii{u_{\delta,h}-u_0}_{L^\infty(\chi_h)}\leq C\delta^2,
\end{equation}
where $C$ is a generic constant independent of $\delta$ {and $h$}.
\end{theorem}

\begin{theorem}[Convergence to the Analytical Nonlocal Diffusion Solution]\label{thm:consistency}
Consider uniform Cartesian grids and $s<d$. Assume that $A(\xb,\yb)\in C^{4}(\overline{(\omg_D)^2})$, the analytical nonlocal diffusion solution $u_\delta(\xb) \in C^1(\overline{\omg_D})$, $a(\xb)\in C^{\infty}(\omg)$ and $\omg_D\in C^1$, then there exists a $\delta_0>0$ and $C_{pos}<1$, such that for a fixed $\delta$ satisfying $0<\delta<\delta_0$ and $h<C_{pos}\delta$, the following convergence property holds for the meshfree quadrature rule with $n=3$:
\begin{equation}
\vertii{u_{\delta,h}-u_\delta}_{L^\infty(\chi_h)}\leq C h^{\min(1,d-s)},
\end{equation}
where $C$ is a generic constant independent of $h$ but may depends on $\delta$.
\end{theorem}


\subsection{Temporal Discretization}

In this section we introduce the discretization methods in time, and demonstrate the fully discretized methods for nonlocal diffusion and peridynamics problems with Dirichlet-type boundary conditions. Although in the demonstrating examples of Section \ref{sec:test} we mostly focus on static nonlocal problems (except the validation problem of the Kalthoff-Winkler experiment), the options and codes of dynamic problems are implemented in OBMeshfree, which can be readily used following the instruction provided in Section \ref{sec:manual}. 

For the dynamic nonlocal diffusion model \eqref{eqn:diffusion_local}, the backward Euler method is employed. With time step size $\Delta t$ and the approximated solution at the $m$-th time step, $u_i^{m}$, the solution at the $m+1$-th time step is solved via:
\begin{equation}\label{eqn:diff_dis}
\left\{\begin{array}{ll}
    \frac{\rho}{\Delta t} u_{i}^{m+1} - (\mcL_{D\delta,h} u)_i^{m+1} =  \fb(\xb_i,(m+1)\Delta t) + \frac{\rho}{\Delta t}u_{i}^m, & \quad \text{for }\xb_i \text{ in }\Omega\bigcap\chi_h,\\
u_i^{m+1}=u_D(\xb_i,(m+1)\Delta t), & \quad \text{for }\xb_i \text{ in }\omgb_D\bigcap\chi_h,\\
u_i^0=\phi(\xb_i), &\quad \text{ for }\xb_i \in\omg_D\bigcap\chi_h,\\
\end{array}\right.
\end{equation}
where $\mcL_{D\delta,h}$ is the discrete nonlocal diffusion operator as defined in \eqref{eqn:discreteNonlocDiff}, $u_D(\xb_i)$ is the prescribed Dirichlet boundary condition, and $\phi$ is the initial value. 

For the dynamic peridynamics model, to discretize in time we also apply the backward time stepping scheme. With time step size $\Delta t$, at the $(m+1)-$th time step we solve for the displacement $\ub_i^{m+1}\approx\ub(\xb_i,(m+1)\Delta t)$ following:
\begin{equation}\label{eqn:probdis}
\left\{\begin{array}{ll}
    \frac{\rho}{\Delta t^2} \ub_{i}^{m+1} - (\mcL_{P\delta,h} \ub)_i^{m+1} =  \fb(\xb_i,(m+1)\Delta t) + \frac{\rho}{\Delta t^2}(2\ub_{i}^m - \ub_i^{m-1}), & \quad \text{for }\xb_i \text{ in }\Omega\bigcap\chi_h,\\
\ub_i^{m+1}=\ub_D(\xb_i,(m+1)\Delta t), & \quad \text{for }\xb_i \text{ in }\omgb_D\bigcap\chi_h,\\
\ub_i^0=\phib(\xb_i),\;\ub_i^1=\psib(\xb_i), &\quad \text{ for }\xb_i \in\omg_D\bigcap\chi_h,\\
\end{array}\right.
\end{equation}
where $\mcL_{P\delta,h}$ is the discretized nonlocal operator as defined in \eqref{eqn:discreteNonlocBond}, $\ub_D$ is the given Dirichlet-type boundary condition, and $\phib$, $\psib$ are the initial displacement field at the first two time steps. 

\subsection{Peridynamics with Free Surfaces and Evolving Fracture}\label{sec:disc_frac}

In this section we extend OBMeshfree, to handle peridynamics models with free surfaces and fracture. For a given point $\xb_i$ and the horizon $\delta$, a bond is associated with each neighbor point $\xb_j\in B_{\delta}(\xb_i)$, and the weight $\omega_{ij}$ is associated with this bond. In the meshfree formulation, the fracture is captured by evolving free surfaces implicitly via the breaking of bonds. When fracture occurs, it creates new surfaces when the free surface boundary conditions are applied. That means, the Neumann-type boundary $\partial\omg_N$ has expended. Instead of parameterizing the $\partial\omg_N$ and evolve its formulation with time, in OBMeshfree the boundary $\partial\omg_N$ is naturally represented by breaking bonds. In particular, when the change of displacement on material point $\xb_j$ may have an impact on the displacement at $\xb_i$, we call their bond as ``intact'', and set the corresponding state function value $\theta(\xb_i,\xb_j,t)$ as $1$. On the other hand, when the bond between $\xb_i$ and $\xb_j$ intersects the surface $\partial\omg_N$, and/or the bond stretch $s(\xb_i,\xb_j,\tau)$ has exceeded the critical bond stretch criteria $s_0(\xb,\yb)$ at some time instant $\tau<t$, the bond is considered ``broken'' and we set $\theta(\xb_i,\xb_j,t)$ as $0$. In particular, at the $m$-th time step we set:
\begin{align}\label{eqn:gamma_disc}
     &\theta^{m}_{ij} = \begin{cases}
    1, \quad \text{if }\xb_j\in B_\delta(\xb_i)\bigcap\Omega_D \text{ and }s(\xb_i,\xb_i,l\Delta t)\leq s_0(\xb_i,\xb_j),\;\forall l=1,\cdots,m,\\
    0, \quad \text{otherwise}, \\
    \end{cases}
\end{align}
Applying the above formulation in \eqref{eq:newform2}, at the $m$-th dynamic step, we seek for solutions of the displacement $u^m_{i}\approx \ub(\xb_i,m\Delta t)$ through the following meshfree scheme:
\begin{align}
   \nonumber &\frac{\rho}{\Delta t^2} \ub_{i}^{m+1}
-  c\sum_{\xb_j \in \chi_h\bigcap B_\delta(\xb_i)\backslash\{\xb_i\}} \kappa(\xb_i,\xb_j) \gamma_\delta(|\xb_i-\xb_j|)\frac{\left(\xb_j-\xb_i\right)\otimes\left(\xb_j-\xb_i\right)}{\left|\xb_j-\xb_i\right|^2}  \left(\ub^{m+1}_i - \ub^{m+1}_j \right) \theta^{m}_{ij}\omega_{ij}\\
&=\fb(\xb_i,(m+1)\Delta t) + \frac{\rho}{\Delta t^2}(2\ub_{i}^m - \ub_i^{m-1}),  \quad \text{for }\xb_i \text{ in }\Omega\bigcap\chi_h,\label{eqn:pd_disc}\\
&\ub_i^{m+1}=\ub_D(\xb_i,(m+1)\Delta t),  \quad \text{for }\xb_i \text{ in }\omgb_D\bigcap\chi_h,\\
&\ub_i^0=\phib(\xb_i),\;\ub_i^1=\psib(\xb_i), \quad \text{ for }\xb_i \in\omg_D\bigcap\chi_h.
\end{align}
Note that because the evolving fracture creates new free surfaces, so $\partial\omg_N$ and $\theta$ alter with $u^{n+1}$. In our implementation, we have been using the damage index, $\theta^{n}_{ij}$, from the last time step, and the above algorithm can therefore be seen as a semi-implicit time integration approach. However, we point out the users can also implement a fully implicit approach  by employing subiterations at each time step following \cite{yu2021asymptotically}, to capture the implicit coupling between the material response and the evolving geometry due to fracture evolution. First one can assume that no new bonds have been broken at the current time step and solve for the displacement field. Second, based on the displacement field, the damage criteria is evaluated and $\theta^{m+1}_{ij}$ is updated following \eqref{eqn:gamma_disc} for each bond. If any bond meets the criteria of breaking, the displacement field will be solved again with new free surfaces. In \cite{yu2021asymptotically}, we repeat this procedure until no new broken bonds are detected, and finally proceed to the next time step.

Finally, the solution $\ub_i^{m+1}$ and the bond state function $\theta^{m+1}_{ij}$ are obtained at time step $m+1$, we postprocess fracture evolution and identify cracks, by calculating the damage field $d_{i}^{m+1}\approx d(\xb_i,(m+1)\Delta t)$ as
\begin{equation}
{d_{i}^{m+1}}=\dfrac{\underset{\xb_j \in \chi_h\bigcap B_\delta(\xb_i)\setminus \xb_i}{\sum}(1-\theta^{m+1}_{ij})}{\underset{\xb_j \in \chi_h\bigcap B_\delta(\xb_i)\setminus\xb_i}{\sum}1},
\end{equation}
which measures the weakening of material via the percentage of broken bonds in the neighborhood of $\xb_i$.

\section{Using OBMeshfree}\label{sec:manual}

In this section we firstly overview the usage of OBMeshfree. For quick start, we refer the readers to Section \ref{sec:usage}. In order to further customize the examples, we then introduce the overall structure of the code and explain each .cpp file in details in Section \ref{sec:structure}. One can find the most up-to-date version of OBMeshfree at \url{https://github.com/youhq34/meshfree_quadrature_nonlocal}.

\subsection{OBMeshfree Usage and Workflow}\label{sec:usage}
Four cases are implement in this exemplar code. After downloading and extract the codes, users can compile the code using the following command

\begin{itemize}
    \item \url{./make} \url{Nldiff} for the static nonlocal diffusion problem, on a default domain $\Omega=[0,1]^2$;
    \item \url{./make} \url{Nldiffd} for the dynamic nonlocal diffusion problem, on a default domain $\Omega=[0,1]^2$;
    \item \url{./make} \url{PD} for the static bond-based peridymics problem, on a default domain $\Omega=[0,1]^2$;
    \item \url{./make} \url{KW} for the dynamic bond-based peridynamics problem with evolving fracture, reproducing the Kalthoff-Winkler fracture experiment.
\end{itemize}
Then, the codes run with:
\begin{itemize}
    \item \url{./nldiff.ex} \url{<#particles>}  \url{<dhratio>} \url{<poly_order>} \url{<case>} for the static nonlocal diffusion problem. Inputs are: 
\begin{itemize}[label=\ding{212}]
    \item \url{<#particles>}: The number of uniform discretization points on each dimension.
    \item \url{<dhratio>}: The ratio between $\delta$ and $h$.
    \item \url{<poly_order>}: The highest polynomial order, $n$, to be exactly reproduced by the quadrature weights\footnote{Here, we point out that the highest polynomial order should be chosen based on the singularity of the kernel. Generally, one should set $n>d+s-3$ according to the truncation estimates.}.
    \item \url{<case>}: A switcher for different experiment settings of Section \ref{sec:test}. Here, \url{0} corresponds to the first example and \url{1} corresponds to the second example. 
\end{itemize}
    \item \url{./nldiffd.ex} \url{<#particles>}  \url{<dhratio>}  \url{<poly_order>} \url{<dt>} \url{<timestep>} for the dynamic nonlocal diffusion problem. The first three inputs are the same as in the static nonlocal diffusion problem. The last two inputs are:
\begin{itemize}[label=\ding{212}]
    \item \url{<dt>}: The time step size.
    \item \url{<timestep>}: The total number of time steps to be simulated.
\end{itemize}
    \item \url{./PMB2D.ex} \url{<#particles>}  \url{<dhratio>}  \url{<poly_order>} \url{<perturbation>} for the static bond-based peridymics problem. The first three inputs are the same as in the nonlocal diffusion problem. The last input parameter is:
\begin{itemize}[label=\ding{212}]
    \item \url{<perturbation>}: The level of perturbation from a uniform grid, so as to create quasi-uniform discretization points. In particular, the Cartesian grids with grid size $h$ are perturbed with a uniformly distributed random vector field  $(\Delta x,\Delta y)$, $\Delta x, \Delta y \sim \mathcal{U}[-rh,rh]$. Here $r$ controls the degree of perturbation, which is given by \url{<perturbation>}.
\end{itemize}
    \item \url{./KW.ex} \url{<#particles>}  \url{<dhratio>}  \url{<poly_order>}
    \url{<dt>} \url{<timestep>} for the dynamic bond-based peridynamics problem with evolving fracture, with the same inputs as in the dynamic nonlocal diffusion problem.
\end{itemize}

\subsection{Software Components}
gcc 7.5.0 or a newer version is required. The codes are built based on BLAS and LAPACK. Users might need to change the settings of BLAS and LAPACK libraries in Makefile. 

\subsection{Structure of the Code}\label{sec:structure}
The current version of OBMeshfree contains three folders, corresponding to nonlocal diffusion problems, bond-based peridynamics problems, and the Kalthoff-Winkler fracture experiment simulation, respectively. 
\begin{itemize}
\item Under the folder \url{nonlocal_diffusion}, there are two .cpp files.  \url{nonlocaldiff_static.cpp} is the script for static nonlocal diffusion problems, 
and \url{nonlocaldiff.cpp} is for dynamic nonlocal diffusion problems.
\item Under the folder \url{bond-based-pd}, \url{PMB_2Dweight.cpp} provides the script for the static bond-based peridynamics model. 
\item Under the folder \url{Kalthoff-Winkler-test}, \url{KW_2Dweight_dynamic.cpp} is the script for the Kalthoff-Winkler fracture experiment, serving as an exemplar script for dynamic fracture problem with peridynamics, as well as a validation of OBMeshfree to realistic engineering applications.
\end{itemize}
In each folder, \url{vvector.h} provides the definitions for basic vector operations.

\subsection{Description of the Code}
In addition to the grid size $h$, the ratio between horizon/grid sizes, $\delta/h$, the reproducing polynomial order $n$, and the time-step size $\Delta t$, users can further customize the test examples and material properties based on their demands. In this section we will illustrate the structure of the scripts and the usage of each function.

\subsubsection{Description of nonlocaldiff\_static.cpp}\label{sec:nonlocal_static_code}
\url{nonlocal_static.cpp} provides numerical simulations for static nonlocal diffusion problems, and serves as a numerical verification for the consistency of numerical solutions to a user-defined analytical local/nonlocal limit. The structures are as follows.

\begin{tcolorbox}[breakable]
\begin{itemize}
\item User defined functions: \begin{itemize}[label=\ding{212}]
\item \url{phi} defines the functions $q$ in $\bm{V}_{h,\xb_i}$, the finite dimensional function space OBMeshfree seeks to exactly reproduced by the quadrature weights.

\item \url{Iphi} defines the analytical integral of each $q$, which will be employed as the right hand side of the constraints in \eqref{eq:quadQP}.

\item \url{inverse} is the function that calculates the inverse matrix, developed based on LAPACK.

\item \url{u_exact} defines the analytical solution at point $\xb=(x,y)$ when it is available, for the purpose of verifying the convergence. This analytical solution will be also used as initial conditions and boundary conditions.

\item \url{Ffun} defines the loading field $f(\xb)$ at point $\xb=(x,y)$, with inputs  $\xb$ and $\delta$.

\item \url{diff_coef} defines the nonlocal diffusion coefficient at point $\xb=(x,y)$ as the harmonic mean of the local diffusion coefficient, when studying the AC convergence. Users can also define the nonlocal diffusion coefficient field directly for a more general nonlocal model, as described below.

\item \url{nonlocal_diff_coef} defines the nonlocal diffusion coefficient $A(\xb,\yb)$, taking a pair of 2D points, $(\xb,\yb)=(x_1,x_2,y_1,y_2)$, as the input.

\end{itemize}
\item \url{Preprocess} performs the main steps of OBMeshfree, i.e. generating the quadrature weights. \url{Preprocess} takes coordinates of the grids, number of grids, $\delta/h$ and the neighborhood list for each grid as inputs, and produces a list of quadrature weights $\omega_{ij}$ corresponding to each neighborhood point for all particles in the domain. 
\item \url{Main} contains the complete procedure solving a static nonlocal diffusion problem. The steps are: 
\begin{enumerate}
\item Read in the number of discretization points in each direction as \url{N} and the ratio $\delta/h$ as \url{dhratio} from inputs.

\item Set initial configuration, including $x$ and $y$ coordinates of the grids, and the analytical solution at each grid.

\item Set neighborhood list for every particle.

\item Generate quadrature weights via \url{Preprocess}.

\item Assemble the stiffness matrix and corresponding forcing term for the particles in the computational domain $\Omega\bigcap\chi_h$.

\item Apply Dirichlet boundary condition for the particles in $\omgb_D\bigcap\chi_h$. 

\item Solve for the linear system using LAPACK build-in functions. 

\item Compute the truncation error and the solution error based on the user-defined analytical solution.
\end{enumerate}
\end{itemize}
\end{tcolorbox}

\subsubsection{Description of nonlocaldiff.cpp}

\url{nonlocaldiff.cpp} runs dynamical simulations for nonlocal diffusion problems. It also provides the option to compare the numerical solutions with a user-defined analytical local/nonlocal limit. The structures are as follows.

\begin{tcolorbox}[breakable]
\begin{itemize}
\item User defined functions:
\begin{itemize}[label=\ding{212}]
\item \url{phi}, \url{Iphi}, \url{diff_coef} and \url{inverse} are defined the same as in Section \ref{sec:nonlocal_static_code}.

\item \url{u_exact} defines the analytical solution at point $\xb=(x,y)$ and time instant $t$.

\item \url{F_fun} defines the forcing term for the nonlocal diffusion equation at point $\xb=(x,y)$ and  time instant $t$.

\end{itemize}
\item \url{Preprocess} generates the quadrature weights and also assembles the stiffness matrix. 

\item \url{Backward_Euler} is the function updating the solution at each time step, with the following steps:
\begin{enumerate} 
\item Update the forcing term and the Dirichlet-type boundary condition for the current time instant $t$. 
\item Based on the solution at the previous time instant, $t-\Delta t$, solve the linear system via LAPACK build-in functions and update the solution at the current time instant $t$.
\end{enumerate}
\item \url{main} contains the complete procedure solving a dynamic nonlocal diffusion problem. The steps are: 
\begin{enumerate}
\item Read in the number of discretization points in each direction as \url{N}, the ratio $\delta/h$ as \url{dhratio}, the time step size $\Delta t$ as \url{dt}, and the total number of time steps as \url{timestep} from inputs.

\item Set up the $x$ and $y$ coordinates of the grids, the initial condition, and the analytical solution at each discretization point.

\item Set up the neighborhood list for every particle.

\item Generate quadrature weights and assemble the stiffness matrix via \url{Preprocess}.

\item Perform the iteration in time from step $1$ to step \url{timestep}. At each time step, run \url{Backward_Euler}.

\item Evaluate the solution error at the last time step.
\end{enumerate}
\end{itemize}
\end{tcolorbox}

\subsubsection{Description of PMB\_2Dweight.cpp}

\url{PMB_2Dweight.cpp} performs the numerical verification for static bond-based peridynamics problem and provides the option to compare the numerical solutions with a user-defined analytical local/nonlocal limit. The structures are as follows.

\begin{tcolorbox}[breakable]
\begin{itemize}
\item User defined functions:
\begin{itemize}[label=\ding{212}]
\item \url{phi} and \url{Iphi} are defined the same as in Section \ref{sec:nonlocal_static_code}.

\item \url{u_exact} and \url{v_exact} define the $x$- and $y$- components of the analytical displacement field at point $\xb=(x,y)$, respectively. 

\item \url{E} defines the heterogeneous material property field, as the function of Young's modulus $E$ at point $\xb=(x,y)$. Then, the nonlocal modulus coefficient is defined as the harmonic mean of the local modulus coefficient.

\item \url{fx_exact} and \url{fy_exact} defines the body force density components of the $x$- and $y$-directions at point $\xb=(x,y)$, respectively. 
\end{itemize}
\item \url{Preprocess} generates the quadrature weights.
\item \url{Main} performs the complete procedure solving a nonlocal static bond-based peridynamics problem. The steps are: 
\begin{enumerate}
\item Read in the number of discretization points in each direction as \url{N} and the ratio $\delta/h$ as \url{dhratio} from inputs.


\item Set up the $x$ and $y$ coordinates of the grids and the analytical solution, if available, at each discretization point.

\item Set up the neighborhood list for every particle.

\item Generate quadrature weights via \url{Preprocess}.

\item Assemble the stiffness matrix and corresponding forcing term for the particles in the computational domain $\Omega \bigcap \chi_h$.

\item Apply Dirichlet boundary condition for the particles in $\omgb_D\bigcap\chi_h$. 

\item Solve for the linear system using LAPACK build-in functions. 

\item Compute the truncation error and the solution error. 
\end{enumerate}
\end{itemize}
\end{tcolorbox}

\subsubsection{Description of KW\_2Dweight\_dynamic.cpp}

\url{KW_2Dweight_dynamic.cpp} performs numerical simulation to reproduce the Kalthoff-Winkler fracture experiment. It provides an exemplar simulation for dynamic bond-based peridynamics problem, and serves as a validation for the applicability of the approach to realistic problems. The structures are as follows.


\begin{tcolorbox}[breakable]

\begin{itemize}
    \item User define functions
    \begin{itemize}[label=\ding{212}]
        \item \url{phi}, \url{Iphi} are defined the same as in Section \ref{sec:nonlocal_static_code}.
        \item \url{BoundaryID} assign each particle an ID, classifies the particles in $\chi_h\bigcap\omgb$ into different regions for further processing. 
        \item \url{Bulk} defines the bulk modulus of the material. 
        \item \url{den} defines the density of the material. 
        \item \url{smax} defines the critical bond stretch of the material. 
        \item \url{u_bc} and \url{v_bc} define the $x$- and $y$-components of the displacement field at point $\xb=(x,y)$ for imposing the boundary conditions in the Kalthoff-Winkler fracture experiment.  \url{u_bc} and \url{v_bc} take the coordinate $\xb=(x,y)$ and the user-defined \url{BoundaryID} of the particle $\xb$ as inputs, to impose the Dirichlet-type boundary condition on the top and the bottom of the object in the Kalthoff-Winkler fracture experiment.
        \item \url{fx} and \url{fy} define the body force density components of the $x$- and $y$-directions at point $\xb=(x,y)$, respectively. 
    \end{itemize}
\item \url{Preprocess} generates the quadrature weights. 
\item \url{Backward_Euler} is the function updating the solution at each time step, with the following steps: 
\begin{enumerate}
            \item Assemble the stiffness matrix using the new quadrature weights and the corresponding forcing term for the particles in the computational domain $\Omega\bigcap\chi_h$. 
            \item Update the boundary condition for $\omgb\bigcap\chi_h$ at time instance $t$. 
            \item Based on the solutions at the previous two time instants, $t-\Delta t$ and $t-2\Delta t$, solve the linear system via LAPACK build-in functions and update the solution at the current time instant $t$.
\end{enumerate}
\item \url{main} contains the complete procedure for simulating the Kalthoff-Winkler fracture experiment. The steps are:
\begin{enumerate}
\item Read in the number of particles in $y$- direction as \url{N}, then set the number of particles in $x$-direction as 2\url{N}. Read in the ratio $\delta/h$ as \url{dhratio}, the time step size as \url{dt}, and the total number of time steps to be simulated as \url{timestep} from input. 

\item Set up material properties, including the bulk modulus \url{Bulk}, critical bond stretch \url{smax}, and material density \url{den} for bond-based peridynamics problems. 

\item Set up the $x$ and $y$ coordinates of the grids and the initial condition at each discretization point.

\item Set up the neighborhood list for every particle.

\item Initialize the pre-notched crack by breaking any bond that intersects with the crack. 

\item Generate quadrature weights through \url{Preprocess}.

\item For $m=1:$ \url{timestep}
\begin{itemize}
\item[7a.] Examine the breaking bonds, update the bond state function values and the quadrature weights. 

\item[7b.] Run \url{Backward_Euler} to update the displacement field.

\item[7c.] Compute the damage field using the current solution

\item[7d.] Save the damage and displacements field for post processing. 
\end{itemize}
\end{enumerate}
\end{itemize}
\end{tcolorbox}

\section{Numerical Examples}\label{sec:test}

In this section, we will use manufactured solutions to test the consistency of OBMeshfree, by investigating the convergence of the numerical solution to the local and nonlocal limits. Four test problems are considered: a static nonlocal diffusion problem verifying the consistency to the analytical nonlocal limit,  a static nonlocal diffusion problem verifying the AC convergence, a static bond-based peridynamics problem verifying the AC convergence, and a dynamic bond-based peridynamics problem reproducing the Kalthoff-Winkler fracture experiment. 
We denote $u_{\delta}$ and $u_0$ as the nonlocal and local analytical solution respectively, and $u_{\delta,h}$ the numerical solution. In example 1, we investigate the convergence of numerical solutions to the nonlocal analytical solution with vanishing $h$ and fixed $\delta$ by calculating $\|u_{\delta,h}- u_\delta \|_{l^2}$, where the root-mean-square takes the form $\vertii{F}_{l^2}:=\sqrt{\dfrac{\sum_{\xb_i \in \chi_h} |F(\xb_i)|^2}{\#(\chi_h)}}$, serving as a numerical approximation of the error in the $L^2(\omg)$ norm. Then in examples 2 and 3, we investigate the convergence of numerical solutions to the local analytical solution under the $\delta$-convergence limit. The differences between the local limit and the nonlocal numerical solution are estimated via $\|u_{\delta,h}- u_0 \|_{l^2}$. Last but not least, in example 4 we consider the Kalthoff–Winkler experiment, wherein the fracture dynamics driven by an impactor striking a pre-notched plate generates an experimentally reproducable crack pattern. These four cases demonstrate the ability of the discretization to resolve both static and dynamic problems.

In nonlocal diffusion examples, we set the kernel $\gamma_\delta$ as a constant without singularity, i.e., $s=0$. In this case, we note due to the $O(\delta^2)$ discrepancy between local and nonlocal diffusion operators, $n=3$ would provide the smallest polynomial space of $\bm{V}_{h,\xb_i}$ to achieve the optimal convergence rate $O(\delta^2)$ in AC tests. Moreover, we further point out that when the grid is uniform, one can equivalently set $n=2$, since the constraint for odd-order polynomials are automatically guaranteed thanks to the symmetric properties. Hence, we set the default reproducing polynomial order for nonlocal diffusion examples as $n=2$. As discussed in Theorems \ref{thm:trunc_err_nonlocal} and \ref{thm:trunc_err}, an $O(\delta^2)$ discrepancy to the local limit under $\delta-$convergence and an $O(h)$ error to the nonlocal limit are anticipated. In peridynamics examples, we employ a singular kernel $\gamma_\delta$ with singularity order $s=1$, and investigate the performance of our solver when the grids are not fully uniform. In this case, we set $n=3$ as the default reproducing polynomial order, since the grids and quadrature weights in each $B_\delta(\xb)$ are no longer symmetric.


\subsection{Examples on Nonlocal Diffusion Problems}

In this section we numerically investigate the convergence properties of OBMeshfree, by studying its performance on two nonlocal diffusion examples with manufactured solutions. To verify the error bounds provided in Theorems \ref{thm:AC} and \ref{thm:consistency}, besides the $l^2$ error we further measure the discrepancy of numerical solution and analytical solution with the $l^\infty$ error given by:
$\vertii{F}_{l^\infty}:=\max_{\xb_i \in \chi_h} |F(\xb_i)|$, serving as a numerical approximation of the error in the $L^\infty(\omg)$ norm.

\subsubsection{Consistency to the nonlocal limit}

\begin{figure}[h!]
\centering
\subfigure{\includegraphics[width=.48\columnwidth]{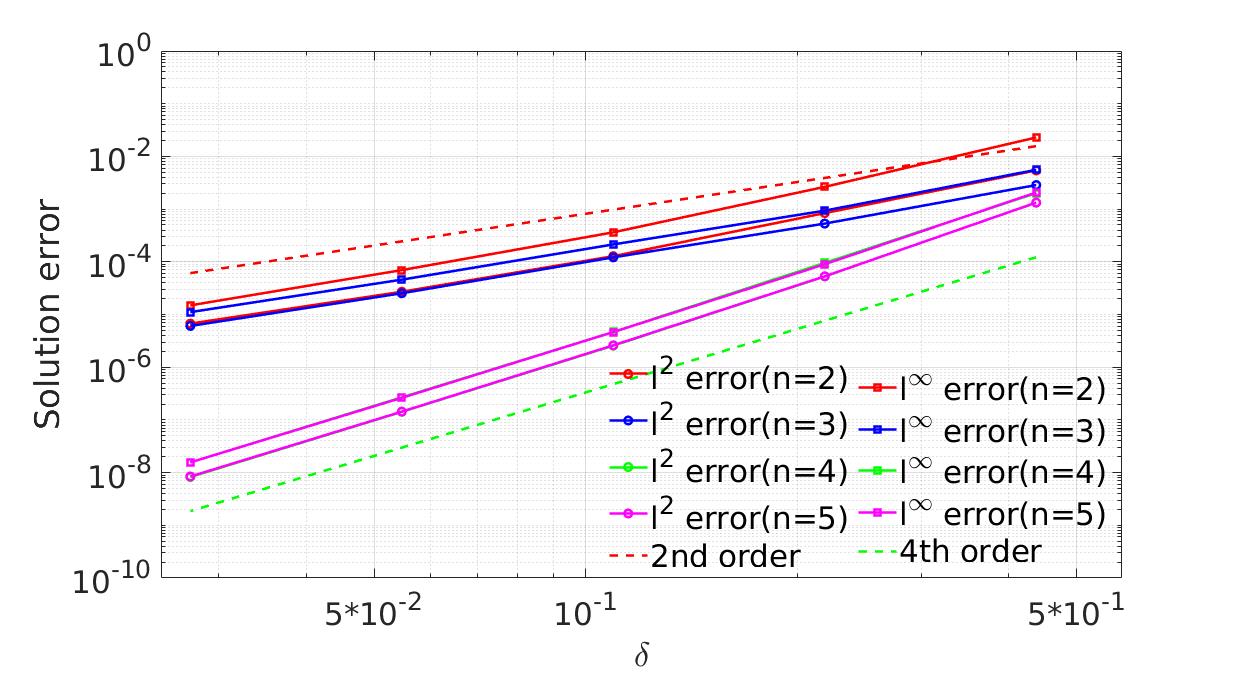}}
\subfigure{\includegraphics[width=.48\columnwidth]{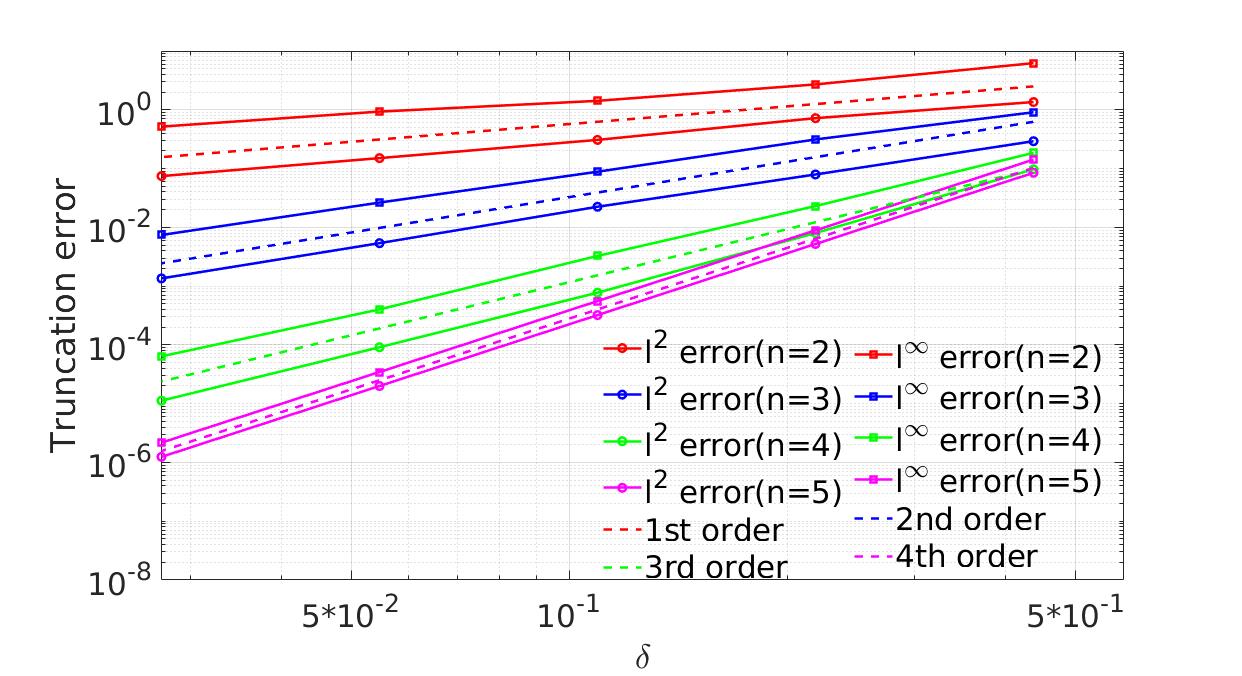}}
\caption{Example 1: verifying the convergence of solution error $\vertii{u_\delta-u_{\delta,h}}$ (left) and truncation error $\vertii{\mcL_{D\delta}[u_\delta]-\mcL_{D\delta,h}[u_\delta]}$ (right) to the \textit{nonlocal limit}, when taking $\delta\rightarrow 0$.
}
\label{fig:diff_nl_ac}
\end{figure}

\begin{figure}[h!]
\centering
\subfigure{\includegraphics[width=.48\columnwidth]{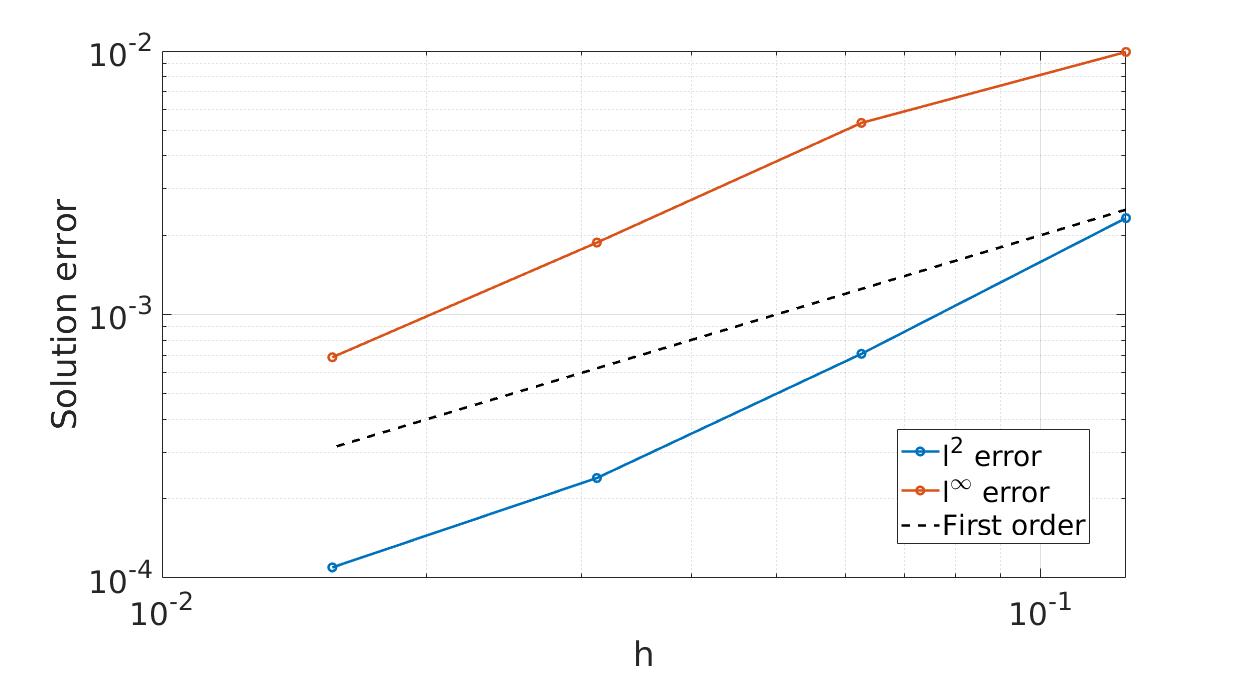}}
\subfigure{\includegraphics[width=.48\columnwidth]{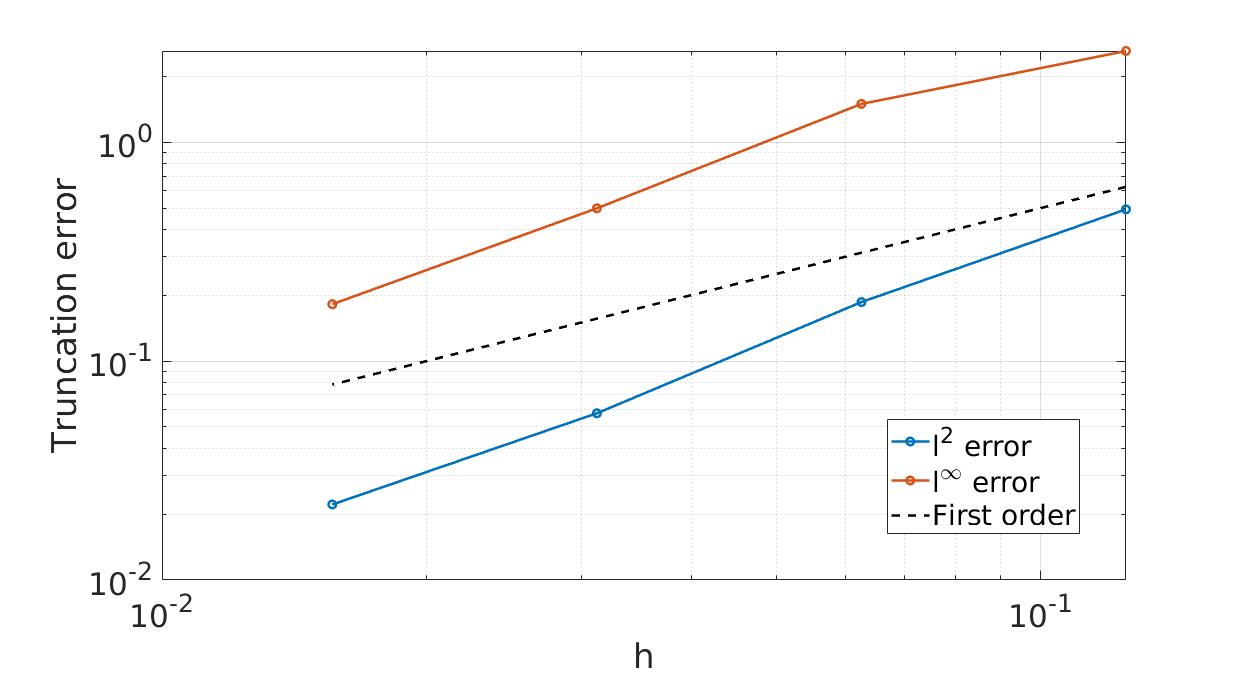}}
\caption{Example 1: verifying the convergence of solution error $\vertii{u_\delta-u_{\delta,h}}$ (left) and truncation error $\vertii{\mcL_{D\delta}[u_\delta]-\mcL_{D\delta,h}[u_\delta]}$ (right) to the \textit{nonlocal limit}, when taking a fixed $\delta$ and setting $h\rightarrow 0$.}
\label{fig:diff_nl}
\end{figure}

Firstly, we test the consistency of numerical solutions to the nonlocal limit. Consider a static nonlocal diffusion problem on $\omg=[0,1]^2$ 
with nonlocal diffusion coefficient field
$$A(\xb,\yb):=A(x_1, x_2,y_1, y_2)=5+x_1+y_1,$$
subjected to the Dirichlet-type boundary condition on $\omgb_D$:
$$u_D(\xb)=u_D(x,y):=x^6+y^6, $$
and a loading field on $\omg$:
$$f(\xb)=f(x,y):=(5+2x)(\frac{5}{32}\pi\delta^8+\frac{15}{8}\pi\delta^6(x^2+y^2)+\frac{15}{4}\pi\delta^4(x^4+y^4))+(\frac{5}{32}\pi\delta^8x+\frac{15}{6}\delta^6x^3+\frac{3}{2}\pi\delta^4x^5).$$ 
This problem has a manufactured analytical nonlocal solution
$$u_\delta(\xb)=u_\delta(x,y)=x^6+y^6.$$

Firstly, we aim to verify Theorem \ref{thm:trunc_err}, by calculating the solution error and truncation error of the discretized nonlocal operator with a fixed $\delta/h=3.5$ and decreasing horizon size $\delta$ from 7/16 to 7/256. The results with reproducing polynomial orders $n=2,3,4,5$ are illustrated in Figure \ref{fig:diff_nl_ac}, where we plot the solution error $\vertii{u_{\delta,h}-u_\delta}$ as well as the truncation error $\vertii{\mcL_{D\delta}[u_\delta]-\mcL_{D\delta,h}[u_\delta]}$ as functions of $\delta$, for each value of $n$. We observe $O(\delta^{n-1})$ convergence in the truncation error, verified the estimate in Theorem \ref{thm:trunc_err}. For the solution error, we observe $O(\delta^n)$ convergence for even $n$ and $O(\delta^{n-1})$ convergence for odd $n$. This is probably caused by the symmetry of the odd-order polynomials on the ball when uniform grids are employed. Users can generate the above results using the script \url{nonlocaldiff_static.cpp}. As an instance, for grid size $h=1/64$, horizon size $\delta=3.5h$, and polynomial order $n=5$, simulations run with the command `\url{./nldiff.ex} \url{64}  \url{3.5}  \url{5} \url{0}', where the last argument corresponds to the case index for this example.

We now proceed to verify Theorem \ref{thm:trunc_err_nonlocal} and Theorem \ref{thm:consistency}, by considering a fixed horizon size $\delta=0.4375$ and decreasing the grid size $h$ from 1/8 to 1/128. When taking the reproducing polynomial order $n=2$ (which is equivalent to $n=3$ due to the symmetry in uniform grids), the results are displayed in Figure \ref{fig:diff_nl}, illustrating an $O(h)$ convergence in both the truncation error and solution error. This result is consistent with the error estimates provided in Theorem \ref{thm:trunc_err_nonlocal} and Theorem \ref{thm:consistency}. Users can generate these numerical results using the script \url{nonlocaldiff_static.cpp}. To keep the horizon size, $\delta$, as a fixed value when decreasing $h$, the second and the third arguments should increase proportionally. For example, to fix the horizon size as $\delta=0.4375$, for grid size $h=1/8$ one should run the command `\url{.\nldiff.ex} \url{8} \url{3.5} \url{3} \url{0}', and for grid size $h=1/16$ one should run `\url{.\nldiff.ex} \url{16} \url{7.0} \url{3} \url{0}'.

\subsubsection{Asymptotic compatibility to the local limit}

\begin{figure}[h!]
\centering
\subfigure{\includegraphics[width=.48\columnwidth]{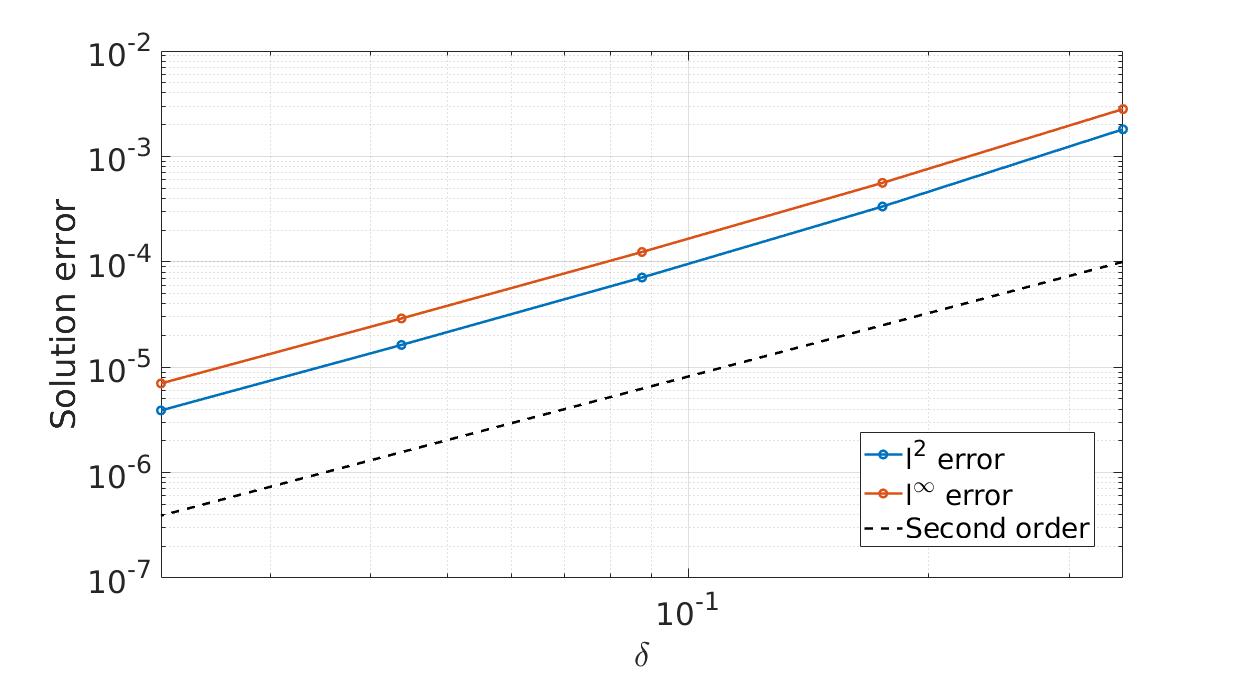}}
\subfigure{\includegraphics[width=.48\columnwidth]{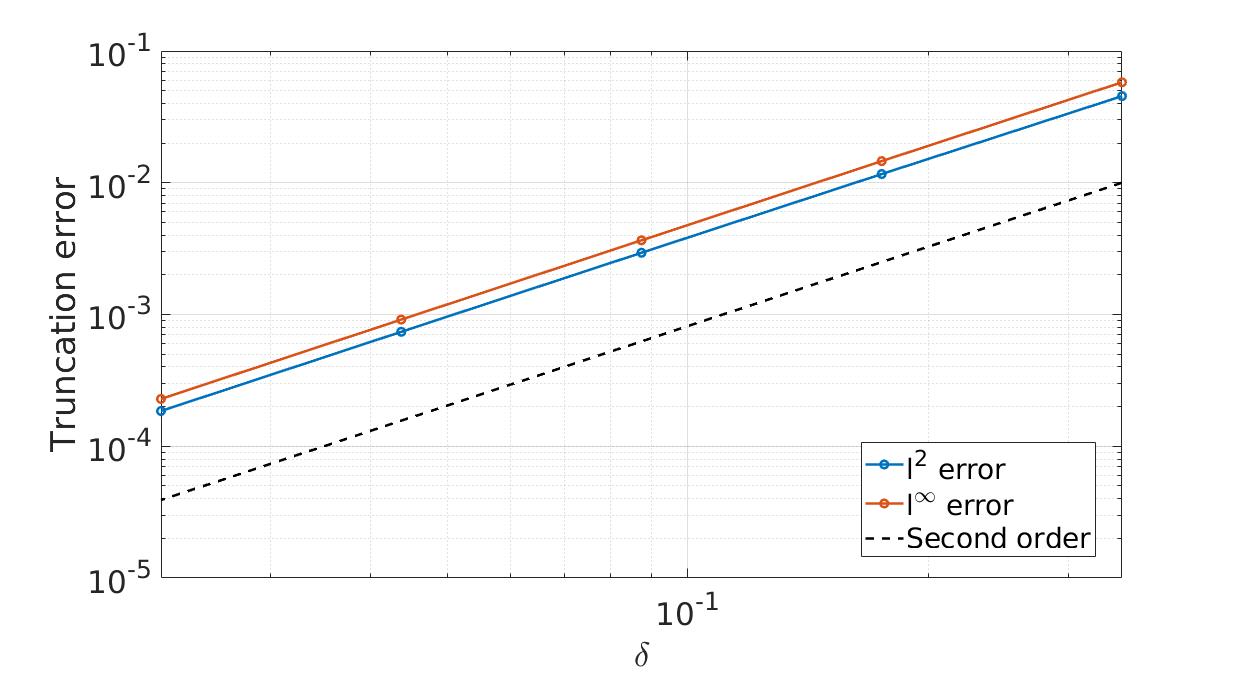}}
\caption{Example 2: verifying the convergence of solution error $\vertii{u_0-u_{\delta,h}}$ (left) and truncation error $\vertii{\mcL_D[u_0]-\mcL_{D\delta,h}[u_0]}$ (right) to the \textit{local limit}, when taking $h,\delta\rightarrow 0$.}
\label{fig:diff}
\end{figure}

Herein, we study the AC consistency of the numerical solution, by considering a static nonlocal diffusion problem example with manufactured local limit. We study a heterogeneous nonlocal diffusion problem on $\omg=[0,1]^2$, with given local diffusion coefficient field 
$$a(\xb)=a(x,y):=2+\sin(x)\sin(y).$$
The object is subjected to Dirichlet-type boundary condition on $\omgb_D$ as:
$$u_D(\xb)=u_D(x,y):=\cos(x)\cos(y)$$
and loading 
$$f(\xb)=f(x,y):=4\cos(x)\cos(y)+4\sin(x)\cos(x)\sin(y)\cos(y).$$
When taking the nonlocal diffusion coefficient field as the harmonic mean of $a(\xb)$ following \eqref{eqn:diffusion_local}, the nonlocal solution should converge to the analytical local solution
$$u_0(\xb)=u_0(x,y)=\cos(x)\cos(y)$$ 
as $\delta\rightarrow 0$.

To present the numerical results of $\delta$-convergence, we fix the ratio $\delta/h=3.5$ and study the convergence of the numerical solution to the local limit with decreasing grid size $h$ from $1/10$ to $1/160$. In this example we employ the reproducing polynomial order $n=2$, which is equivalent to $n=3$ due to the symmetry in uniform grids, and therefore is anticipated to provide the optimal convergence rate, $O(\delta^2)$, to the local limit. The solution error $\vertii{u_{\delta,h}-u_0}$ as well as the truncation error $\vertii{\mcL_D[u_0]-\mcL_{D\delta,h}[u_0]}$ are plotted as functions of $\delta$ in Figure \ref{fig:diff}. In both $l^2$ and $l^\infty$ norms, a $O(\delta^2)$ convergence is observed for the solution error and the truncation error, verifying the analysis in Theorem \ref{thm:AC}. To reproduce these results, users can run the script \url{nonlocaldiff_static.cpp}. Taking the case with grid size $h=1/160$, reproducing polynomial order $n=2$, and $\delta/h=3.5$ as an instance, results are obtained with the command `\url{./nldiff.ex} \url{160} \url{3.5} \url{2} \url{1}'. Here, the last argument corresponds to the case index for this example.


\subsection{Examples on Peridynamics}

In this section we demonstrate two bond-based peridynamics examples using OBMeshfree. In the first example, we provide verification on a static example with manufactured local solution. To provide heuristic studies on the solution convergence rates, we measure discrepancy of numerical solution and analytical solution with the $l^2$ error. Then, in the second example we use the dynamic peridynamics code to reproduce the Kalthoff–Winkler fracture experiment, demonstrating the applicability of this approach to realistic engineering applications involving dynamic fracture.

\subsubsection{Asymptotic compatibility to the local limit}

\begin{figure}[h!]
\centering
\subfigure{\includegraphics[width=.48\columnwidth]{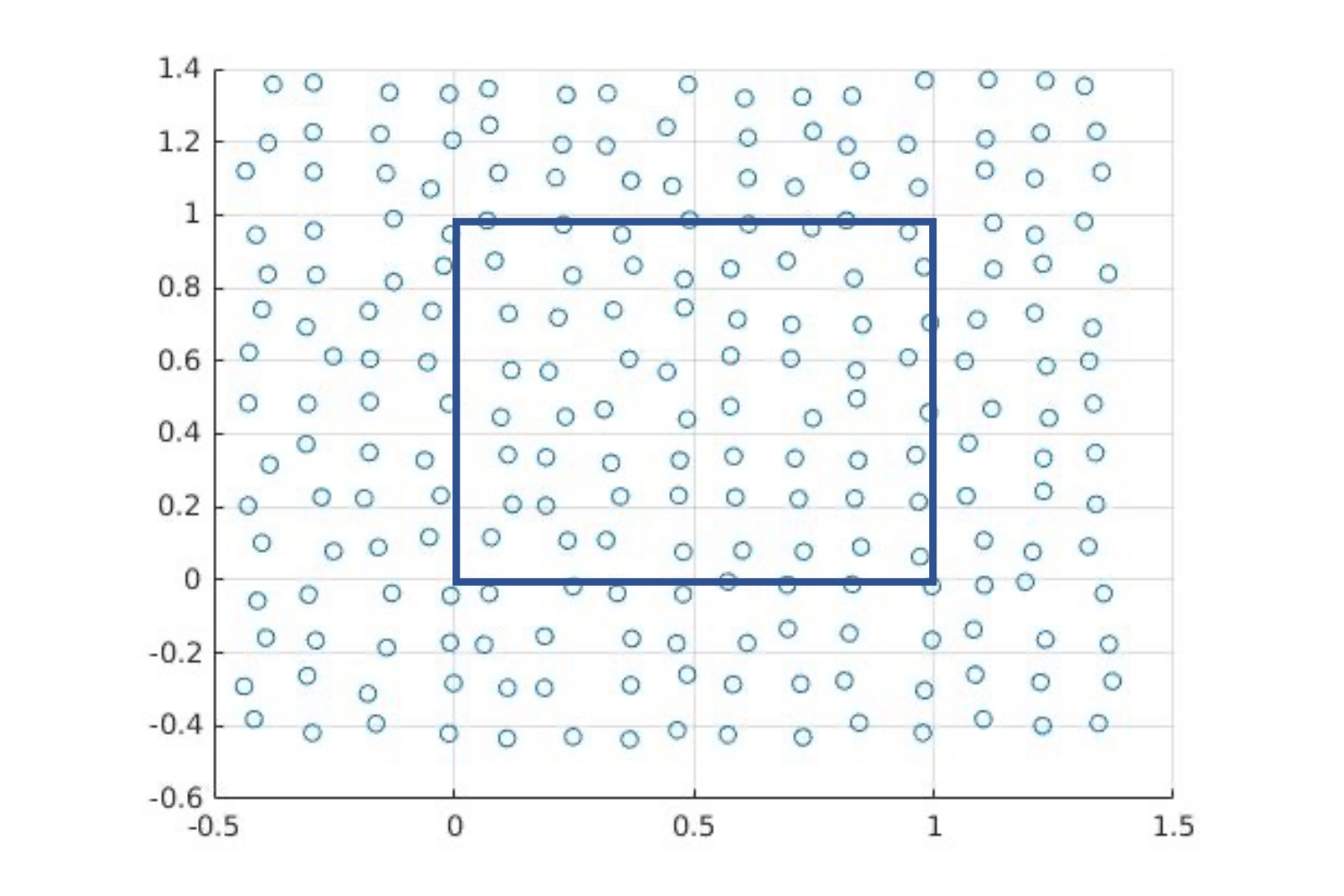}}
\subfigure{\includegraphics[width=.48\columnwidth]{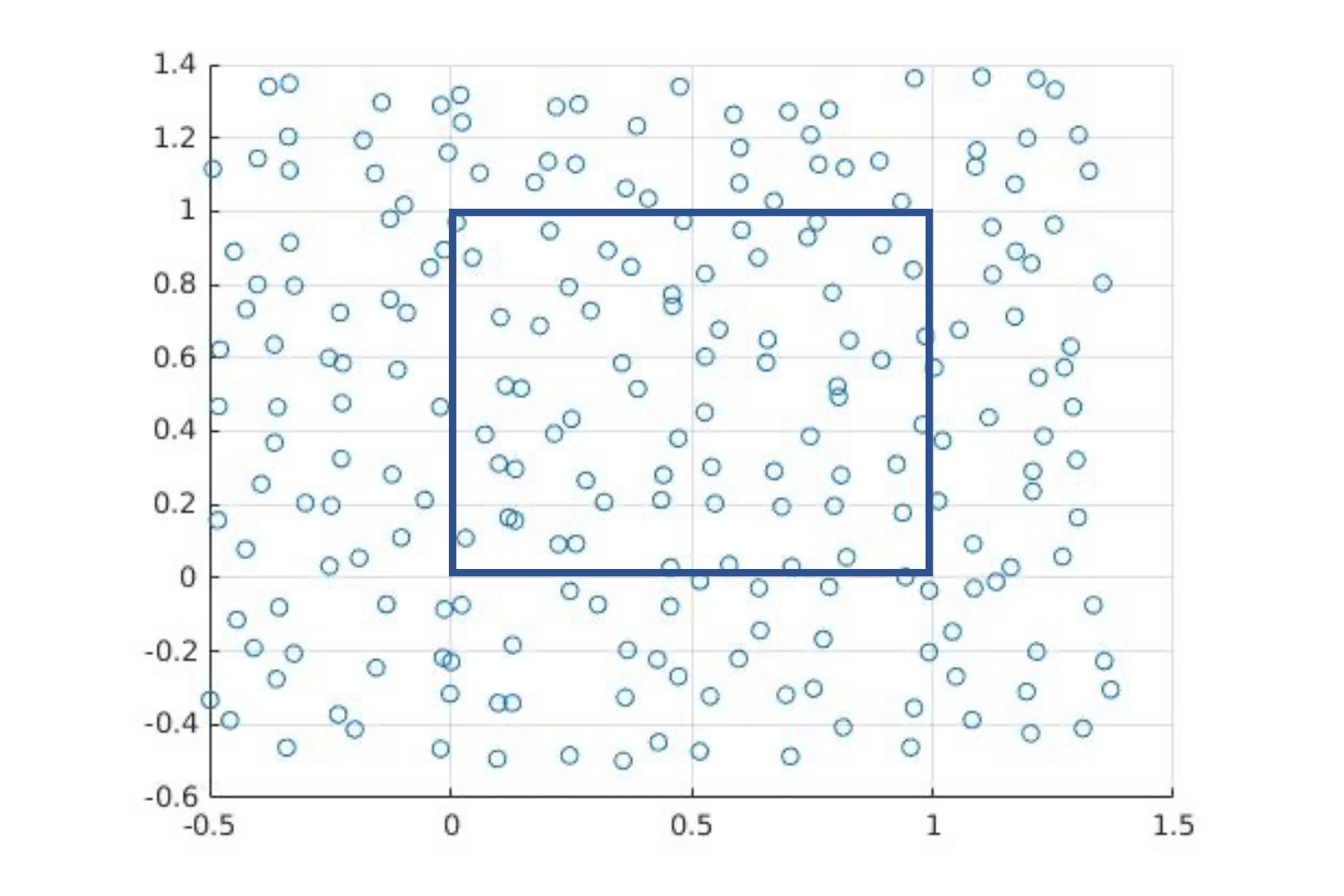}}
\caption{Exemplar non-uniform grids generated for $\omg\bigcup\omgb$, with the perturbation ratios $r=20\%$ (left) and $r=50\%$ (right). The computational domain $\omg=[0,1]^2$ is indicated by the blue box. }
\label{fig:peri_grid}
\end{figure}

\begin{figure}[h!]
\centering
\subfigure{\includegraphics[width=.48\columnwidth]{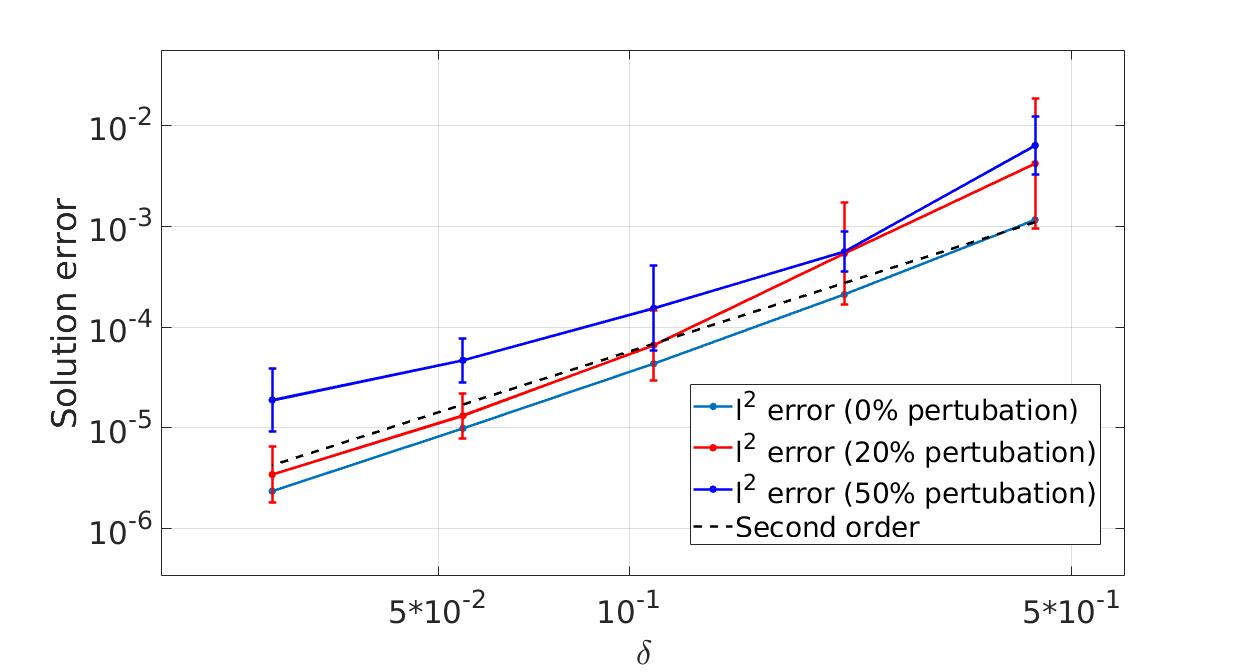}}
\subfigure{\includegraphics[width=.48\columnwidth]{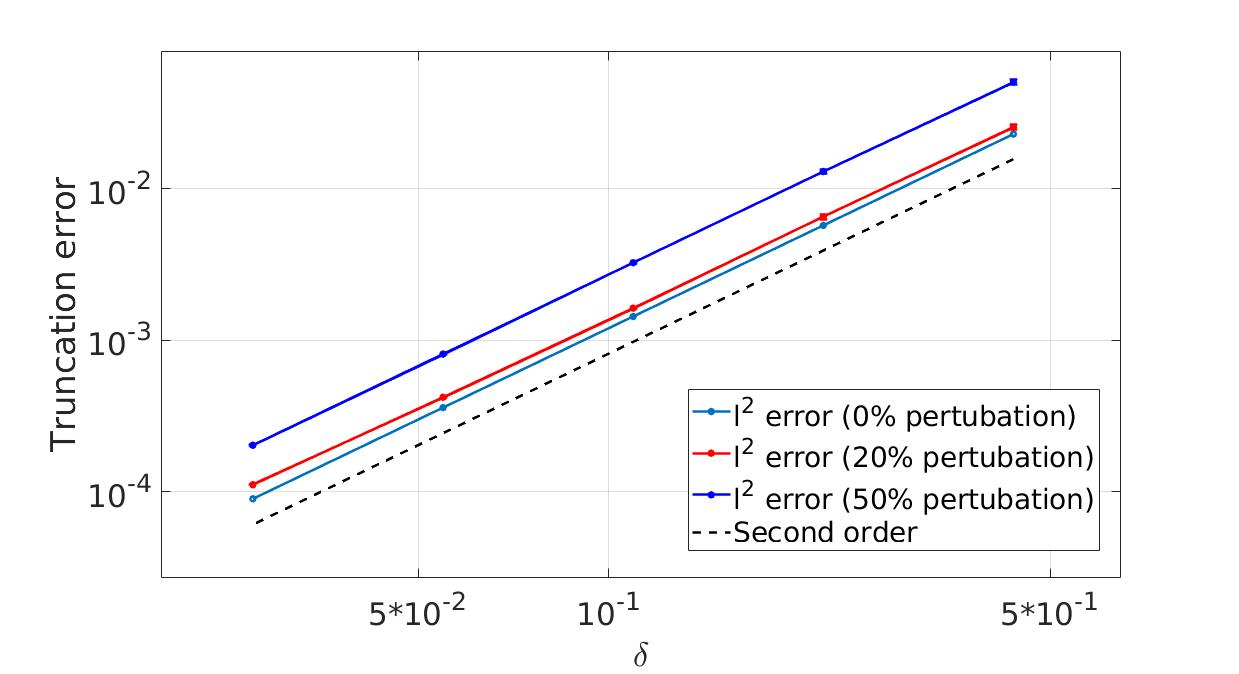}}
\caption{
Example 3: verifying the convergence of peridynamics solution error $\vertii{u_0-u_{\delta,h}}$ (left) and truncation error $\vertii{\mcL_P[u_0]-\mcL_{P\delta,h}[u_0]}$ (right) to the \textit{local limit}, when taking $h,\delta\rightarrow 0$ and perturbing the uniform grid with different levels of perturbation. Means and standard errors of 5 realizations are reported for each perturbation level.
}
\label{fig:peri}
\end{figure}


In this section we consider the static bond-based peridynamics modeling for an object occupying the region $\omg=[0,1]^2$, whose material microstructure is characterized by a fixed Possion ratio $\nu=0.25$ and a local Young's modulus field
$$E(\xb)=E(x,y):=2+\sin(x)\sin(y).$$
Assume that the Dirichlet-type boundary condition
$$\ub_D(\xb)=\ub_D(x,y):=\left[\sin(x)\sin(y),-\cos(x)\cos(y)\right]$$
is given for $\xb\in\omgb_D$, and the object is subject to a body load
\begin{equation*}
\fb(\xb)=\fb(x,y):= \left[\begin{array}{c}
-12C_1\sin(x)\sin(y)+4C_1\cos(2x)\sin^2(y))+4C_1\cos(2y)\sin^2(x))\\
12C_1\cos(x)\cos(y)+3C_1\sin(2x)\sin(2y)
\end{array}\right]^T, \; C_1:=1/(2(1+\nu)).
\end{equation*}
When taking the nonlocal modulus field as the harmonic mean of the local one, the bond-based peridynamics problem converges to the Navier equations \cite{mengesha2012nonlocal,mengesha2014bond,mengesha2014nonlocal} for linear elasticity:
$$-\mcL_{P}[\ub](\xb)=\fb(\xb)$$
as $\delta\rightarrow 0$. Here, 
\begin{equation}\label{eqn:local}
\mcL_{P} [\ub]:=\dfrac{E}{2(1+\nu)} \nabla\cdot(2\mathbf{E}+\text{tr}(\mathbf{E})\mathbf{I}),
\end{equation}
where the strain tensor $\mathbf{E}:=\dfrac{1}{2}(\nabla \ub+(\nabla \ub)^T)$ and $\text{tr}(\mathbf{E})$ denotes its trace. Hence, with this example we aim to investigate if the numerical solution would converge to the following analytical local solution of \eqref{eqn:local}:
$$\ub_0(\xb)=\ub_0(x,y):=\left[\sin(x)\sin(y),-\cos(x)\cos(y)\right]$$
under the $\delta$-convergence setting.

To study the AC convergence, in this example we fix $\delta/h=3.5$ while decrease the grid size $h$ from 1/8 to 1/128, and calculate the discrepancy between our numerical solution and the analytical local limit. To demonstrate our meshfree approach in handling non-uniform grids, we first generate the uniform grid, then perturb the uniform grid points with $(\Delta x,\Delta y)$, where $\Delta x, \Delta y \sim \mathcal{U}[-rh,rh]$. Here, the perturbation ratio $r\in(0,1)$ provides a metric for the effect of anisotropy in the underlying discretization. Two examples of perturbed grids are demonstrated in Figure \ref{fig:peri_grid}, corresponding to $r=0.2$ and $r=0.5$, respectively. For each perturbation ratio $r\in \{0.2,0.5\}$, we generate 5 realizations of grids using different random seeds. To investigate the impact of non-uniform grids, we record the solution and truncation errors for each realization, and report their means and standard errors versus the horizon size $\delta$ in Figure \ref{fig:peri}.  $O(\delta^2)$ convergence is observed in the truncation error for all perturbation levels. Here we notice that the standard errors of the right plot are almost invisible, showing that the truncation errors are not sensitive to perturbations in the discretization grids, possibly because they are dominated by the discrepancy from $\ub$ to the reproducing polynomial space instead of the interpolation error. For the solution error, we also observe $O(\delta^2)$ convergence to the local limit when using uniform grids, while the convergence rate slightly deteriorates as we increase the level of perturbation, possibly due to the effects of grid anisotropy on the stiffness matrix. To reproduce these results, users can run the script \url{PMB_2Dweight.cpp}. For grid size $h=1/128$, $\delta/h=3.5$, and perturbation level $r=20\%$, as an instance, results are generated using the command `\url{./PMB2D.ex} \url{128} \url{3.5} \url{3} \url{0.2}'.

\subsubsection{Dynamic fracture: reproducing the Kalthoff-Winkler experiment}\label{sec:frac}

\begin{figure}[htp]
    \centering
    \includegraphics[width=0.9\columnwidth]{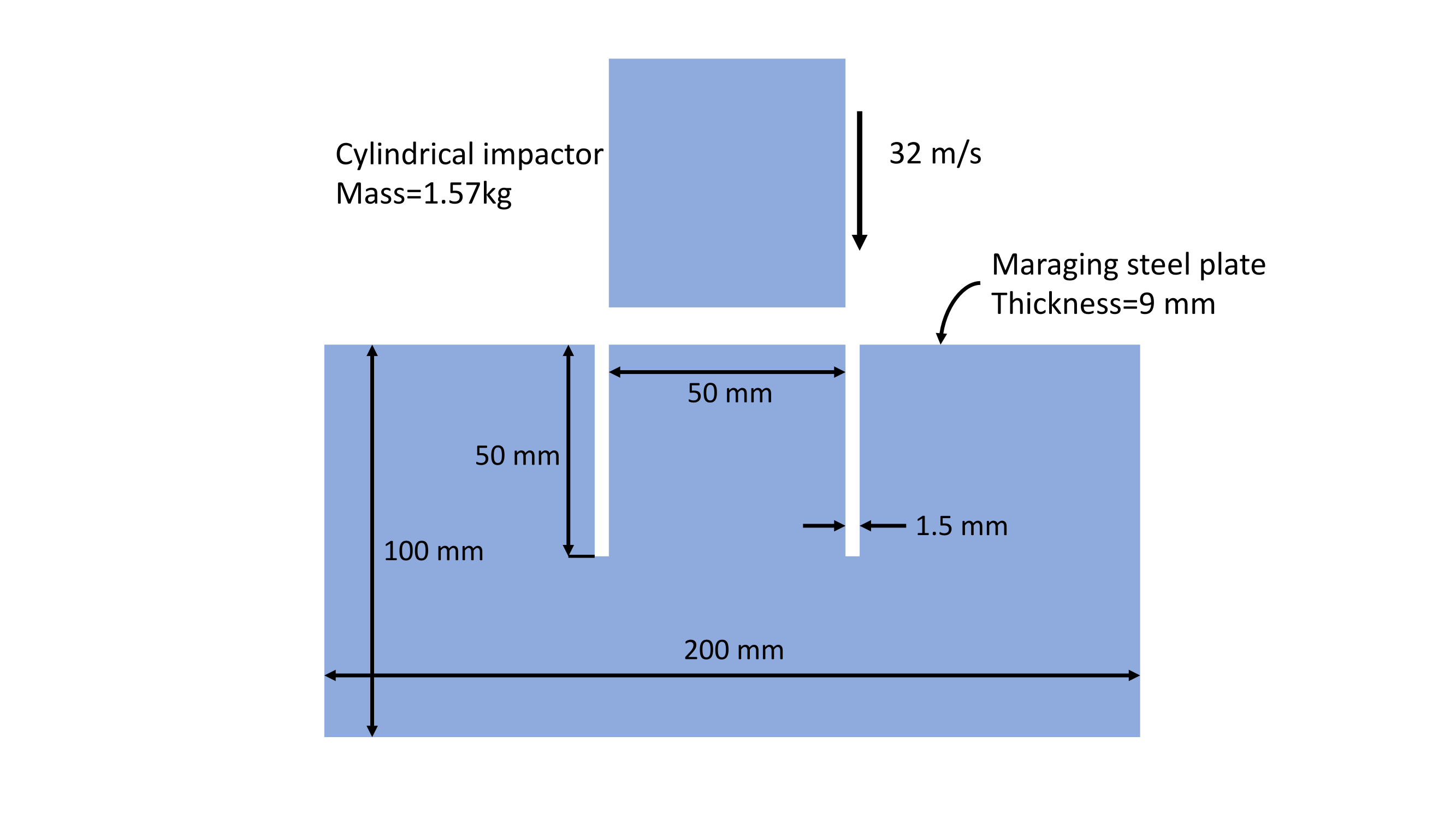}
    \caption{Experimental setup for the Kalthoff-Winkler experiment. The figure is adapted from \cite{silling2001kalthoff}.}
    \label{fig:KWsetup}
\end{figure}


We now consider a dynamic fracture problem where a steel plate is struck by a cylindrical impactor \cite{kalthoff1988failure}, which is the so-called  Kalthoff-Winkler experiment. The plate is pre-notched, and crack will grow from the pre-notch tips upon an impact. Experiments show that the fracture pattern behaves differently depending on the regimes governed by the impactor velocity. In this example, we employ parameters given in \ref{fig:KWsetup}, which is also consistent with those investigated previously in a particle-based peridynamics model \cite{silling2001kalthoff}. Under such a setting,  experimentally it was observed that a reproducible 68$\degree$ angel is formed by the growing crack and the initial vertical pre-notch \cite{kalthoff1988failure}. With this example, we aim to validate if our OBMeshfree is capable to capture the evolving fracture and reproducing the crack pattern.

In this example, to model the impact, $\ub=[0,-32t]$ is imposed between the two notches, as depicted in \ref{fig:KWsetup}. Then, on the rest regime of the top of the plate, a homogeneous Dirichlet boundary conditions $\ub=[0,0]$ is applied. All other boundaries are treated as free surfaces, hence any bonds across those surfaces are set as broken following Section \ref{sec:frac}. For the material properties, we employ the settings employed in \cite{silling2001kalthoff}. In particular, the plate has a density 8e-3 kg/m$^3$, an elastic modulus of 191 GPa, a yield strength of 2000 MPa, and a fracture toughness of 90 MPa m$^{1/2}$. The above material properties yield a bond breaking criteria of $s_0=0.0099/\sqrt{\delta}$ following \eqref{eq:criteria}. In our code, the spatial unit is unified as cm, the temporal unit is unified as ms, and weight unit is unified as kg.

\begin{figure}[h!]
\centering
\subfigure[t=2e-4ms]{\includegraphics[width=.40\columnwidth]{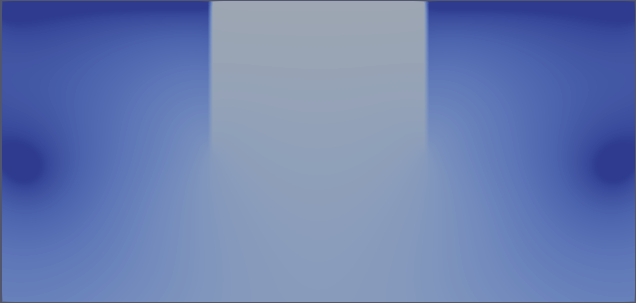}}
\subfigure[t=2e-3ms]{\includegraphics[width=.40\columnwidth]{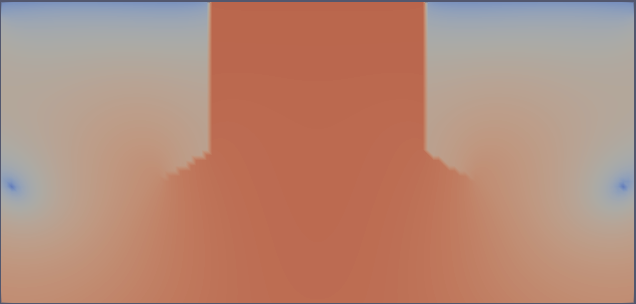}}
\subfigure[t=4e-3ms]{\includegraphics[width=.40\columnwidth]{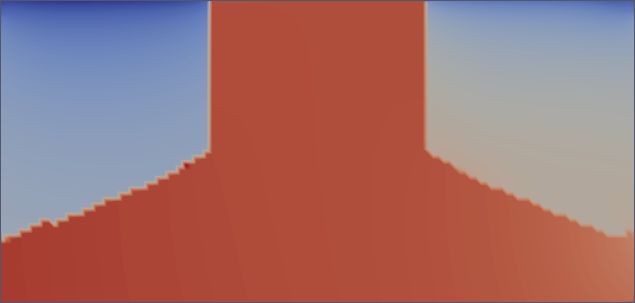}}
\subfigure[t=6e-3ms]{\includegraphics[width=.40\columnwidth]{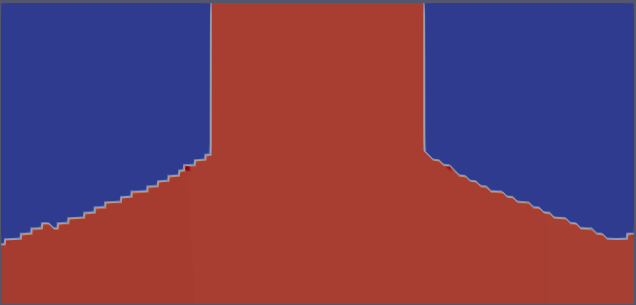}}
\caption{The evolution of the displacement field for the Kalthoff-Winkler fracture experiment, after 1, 50, 100 and 150 time steps. The logarithmic value of displacement magnitude ($\log_{10}\verti{\ub}$) is colored for plot.}
\label{fig:KWdis}
\end{figure}

\begin{figure}[htp]
    \centering
    \includegraphics[width=0.75\columnwidth]{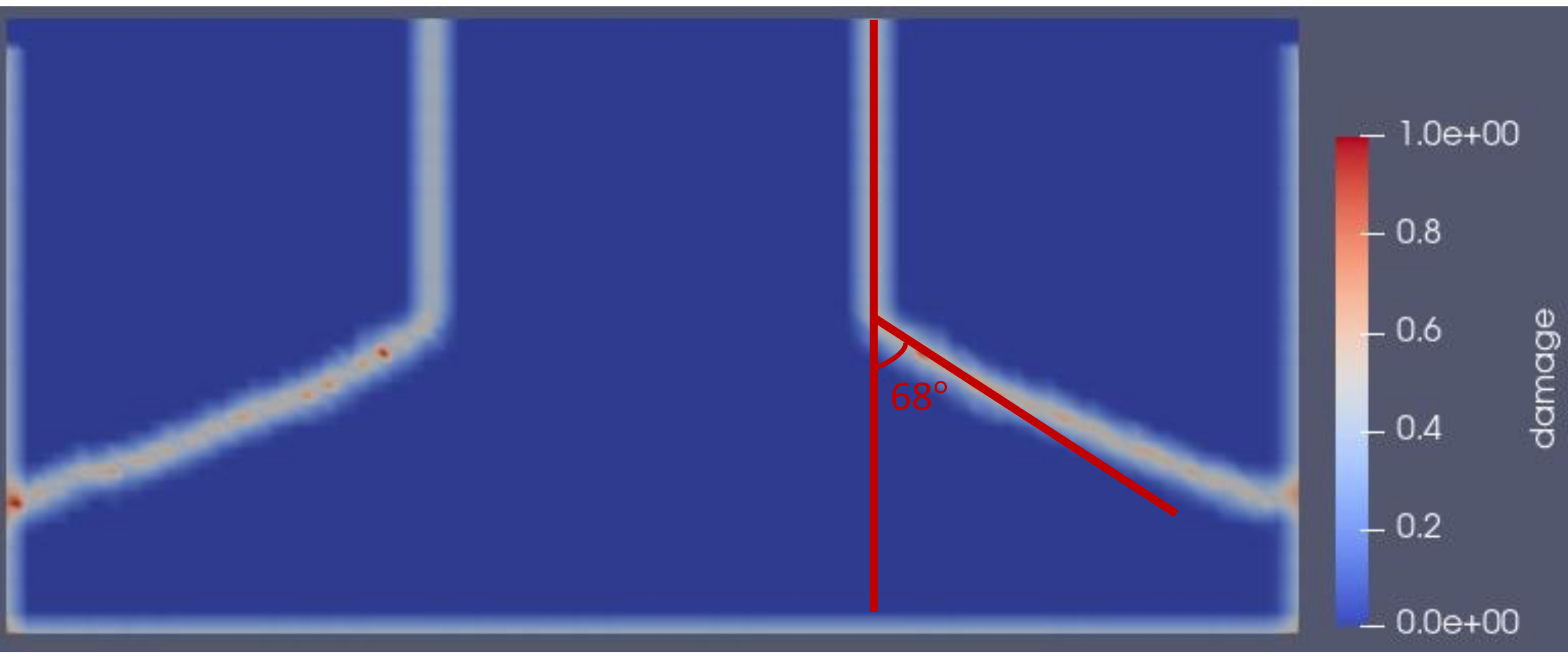}
    \caption{Crack pattern for Kalthoff-Winkler experiment when taking $h=$0.15625cm and $\delta=3.0h$ at time t=6e-3ms, successfully reproduced the 68$\degree$ crack angle as reported in  \cite{kalthoff1988failure,silling2001kalthoff}.}
    \label{fig:KWcrack}
\end{figure}

In Figure \ref{fig:KWdis} we illustrate the evolution of simulated displacement field at four time instances: t=2e-4 ms, t=2e-3 ms, t=4e-3 ms, and t=6e-3 ms. In this simulation a uniform grid with $64\times 128$ particles, $\delta=3.0h$ and time discretization size $\Delta t=2e-4ms$ are employed. At the end of the simulation, three fragments remain due to the crack. In Figure \ref{fig:KWcrack} we further illustrate the fracture pattern,  whose fragment shape reproduces the experimentally observed 68$\degree$ crack angle at the pre-notch tip. The results in this section are generated based on the script \url{KW_2Dweight_dynamic.cpp}. Users can reproduce these results by running the command `\url{./KW.ex} \url{64} \url{3.0} \url{3} \url{2e-4} \url{500}'.


\noindent\textbf{Remark}: To fully resolve the transient dynamics of the problem in \cite{ren2016dual}, we may define the $CFL$ condition number $C_{CFL}:=\dfrac{C_R \Delta t}{h}$, where $C_R$ is the Rayleigh speed calculated following \cite{graff2012wave}. When $C_{CFL}\leq\frac{1}{2}$, we fully resolve the transient dynamics of the problem. In this example, our settings are corresponding to $C_{CFL}\approx0.4$. For further results with different $CFL$ condition numbers, we refer interested readers to \cite{trask2019asymptotically}, where both fully resolved dynamics ($C_{CFL}\leq\frac{1}{2}$) and implicit solution of wave propagation ($C_{CFL}>1$) were investigated.

\section{Conclusion}\label{sec:conclusion}

In this work, we have developed an open-source software called OBMeshfree for meshfree analysis on nonlocal problems. The program is developed based on an optimization-based quadrature rule and consists of a set of routines for generating quadrature weights on a neighborhood of each material point, discretizing two-dimensional nonlocal diffusion and peridynamics operators, performing integrating in time, handling material heterogeneity and evolving fracture, and calculating solution and truncation errors in cases with manufactured solutions. Benchmark problems are presented to verify the convergence, efficiency, and robustness properties of the meshfree discretization method  implemented in OBMeshfree, under both uniform and highly non-uniform nodal distributions. Our method features a unified mathematical workflow for handling material heterogeneity and evolving material fracture, and it is able to provably obtain high-order convergence to both local and nonlocal limits. With sufficient regularity assumptions on the solution and material property fields, the approach is able to obtain $O(h)$ and $O(\delta^2)$ convergences to the nonlocal and local theory, respectively.

The open source code can serve as an entry point for researchers who are interested in the computer implementation of the optimization-based quadrature rule employed in \cite{d2022prescription,fan2021asymptotically,fan2022meshfree,foss2022convergence,trask2019asymptotically}, and the code can also be adopted as a rapid prototyping and testing tool for the simulations with nonlocal models. Although the linear diffusion and bond-based peridynamics problems are chosen as the model problem, the flexibility of the code allows the extension to solve more advanced nonlocal problems for various scientific and engineering applications.

\section*{Conflict of Interest Statement} 
 
On behalf of all authors, the corresponding author states that there is no conflict of interest.

\end{document}